\documentclass[10pt]{amsart}
\usepackage{graphicx}
\usepackage{cite,amssymb,amsmath}
\usepackage{amsthm}

\theoremstyle{plain}

\newtheorem{theorem}{Theorem}[section]
\newtheorem{proposition}[theorem]{Proposition}
\newtheorem{corollary}[theorem]{Corollary}

\newtheorem*{maintheorem}{Theorem~\ref{mapping torus hyperbolic}}

\theoremstyle{remark}
\newtheorem*{remark}{Remark}
\newtheorem*{remarks}{Remarks}

\newcommand{\newoperator}[2]{\DeclareMathOperator{#1}{#2}}
\renewcommand{\H}{\qopname\relax o{H}}
\newoperator{\AH}{AH}
\newoperator{\AJ}{AJ}
\newoperator{\BS}{BS}               
\newoperator{\GF}{GF}               
\newoperator{\GH}{GH}
\newoperator{\GJ}{GJ}
\newoperator{\Isom}{Isom}
\newoperator{\KS}{KS}               
\newoperator{\ML}{ML}
\newoperator{\Nbhd}{Nbhd}
\newoperator{\PL}{PL}
\newoperator{\QF}{QF}
\newoperator{\QH}{QH}
\newoperator{\QJ}{QJ}
\newoperator{\QpH}{Q^\prime H}
\newoperator{\Rad}{Rad}
\newoperator{\Tr}{Tr}
\newoperator{\Vol}{Vol}
\newoperator{\area}{area}
\newoperator{\base}{base}
\newoperator{\codim}{codim}
\newoperator{\const}{const}
\newoperator{\degree}{degree}       
\newoperator{\glue}{glue}
\newoperator{\graft}{graft}
\newoperator{\inj}{inj}             
\newoperator{\interior}{int}
\newoperator{\length}{length}
\newoperator{\llength}{\underline {length}}
\newoperator{\ls}{ls}
\newoperator{\qf}{qf}
\newoperator{\skin}{skin}
\newoperator{\teich}{\textsl{T}}
\newoperator{\unglue}{unglue}
\newoperator{\wb}{wb}               
\newoperator{\window}{window}
\newoperator{\ws}{ws}               

\newcommand{\hy}{\mathbb{H}}
\newcommand{\integers}{\mathbb{Z}}
\newcommand{\complexes}{\mathbb{C}}
\newcommand{\proj}{\mathbb{P}}
\newcommand{\reals}{\mathbb{R}}

\newcommand{\refin}[2]{#1 of \cite{#2}}
\newcommand{\fund}[1]{\pi_1(#1)}
\newcommand{\set}[1]{\{#1\}}

\newcommand{\arrow}{\rightarrow}
\newcommand{\boundary}{\partial}
\newcommand{\cinfty}{S_\infty^1}
\newcommand{\compose}{\circ}
\newcommand{\cross}{\times}

\newcommand{\degrees}{^\circ}
\newcommand{\homeomorphic}{\approx}
\newcommand{\homotopic}{\simeq}
\newcommand{\hs}{\AH(\bar S \cross I, \boundary \bar S \cross I)}

\newcommand{\intersect}{\cap}
\newcommand{\inverse}{^{-1}}

\newcommand{\sinfty}{S_\infty^2}
\newcommand{\thin}{{\hbox{\rm thin}}}
\newcommand{\union}{\cup}
\newcommand{\zeemod}{\integers_}

\begin{document}


\title[Hyperbolic Structures II: Fibered case]{Hyperbolic Structures on 3-manifolds, II:  \\
Surface groups and 3-manifolds which fiber over the circle}
\author{William P. Thurston}
\address{
Mathematics Department \\
University of California at Davis \\
Davis, CA 95616
}
\email{wpt@math.ucdavis.edu}
\date{August 1986 preprint $\to$ Jan 1998 eprint}
\thanks{This project has been supported by the NSF,
currently \#DMS-9704135}
\begin{abstract}
The main result (\ref{mapping torus hyperbolic}) of this paper
is that every atoroidal three-manifold that fibers over the
circle has a hyperbolic structure. Consequently,
every fibered three-manifold admits a geometric decomposition.
The main tool for constructing hyperbolic structures on fibered
three-manifolds is the double limit
theorem (\ref{double limit theorem}), which is of interest for its own sake
and lays out general conditions under which sequences of quasi-Fuchsian
groups have algebraically convergent subsequences.
The main tool in proving the double limit theorem
is an analysis of the geometry of hyperbolic manifolds that are
homotopy equivalent to a surface.
This analysis is also of interest in its own right.

This eprint is based on the August 1986 version of this preprint, which
was submitted, refereed, and accepted for publication; for reasons
that are hard to fathom, I never returned a corrected version to the journal. 
I apologize for my long neglect of its publication, and I want to
thank the referee for detailed comments
which have been incorporated into the present eprint. No other
significant changes have been made, except conversion to \LaTeX,
which has resulted in changes in numbering.
The 1986 preprint was in turn a revision of a 1981 preprint; the
various versions were fairly widely circulated in the early 1980's,
and the results became widely known and used.

Too many developments have intervened 
to be easily summarized, except for pointers particularly
to the works of McMullen
\cite{McMullen:fibered} and Otal \cite{Otal:fibered} that
give alternative proofs for the main results of this paper
and contain other interesting material as well.
\end{abstract}
\maketitle
\begin{figure}[hb]
\centering
\includegraphics{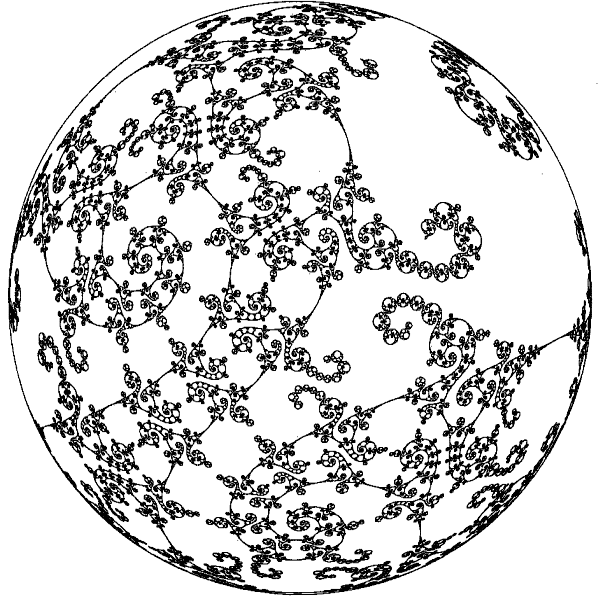}
\caption{This is a portion of the limiting sphere-filling curve
for a fiber of the punctured torus bundle over $S^1$ with gluing
map $R^4 L$, where $R$ is a right-handed Dehn twist about a
$(1,0)$-curve and $L$ is a left-handed Dehn twist around a $(0,1)$-curve.
}
\label{4-1 limit}
\end{figure}
\tableofcontents

\setcounter{section}{-1}
\section{Introduction}

This is the second in a series of papers dealing with the conjecture that all
compact 3-manifolds admit canonical decompositions into geometric pieces.
The main purpose of the current paper is to prove

\begin{theorem}[Mapping torus hyperbolic]
\label{mapping torus hyperbolic}
Let $M^3$ be a compact $3$-manifold (possibly with boundary) which
fibers over
$S^1$, and whose fiber is a compact surface of negative Euler characteristic.

Then the interior of $M$ either
\begin{description}
\item[(i)] has a complete $\hy^2 \cross \reals$ structure of finite volume,
and can be described as a Seifert fibration over some hyperbolic $2$-orbifold,
\item[(ii)] contains an embedded incompressible
torus not isotopic to a boundary component,
and splits along this torus into two simpler three-manifolds, or
\item[(iii)] (generic case) has a complete hyperbolic structure of finite volume.
\end{description}
Cases (i) and (ii) are not mutually exclusive, but (iii) excludes
the other two cases.
\end{theorem}

If $M$ fibers over the circle and has non-empty boundary, then the boundary is
a union of tori and (in the non-orientable case) Klein bottles.

A statement equivalent to the main theorem is that a $3$-manifold admits a
hyperbolic structure if and only if it is homotopically atoroidal. A
$3$-manifold is {\it homotopically atoroidal} if every map of a torus into the
manifold which is injective on fundamental groups is homotopic to the
boundary.  This is not quite the same as the condition that every embedded
incompressible torus is isotopic to the boundary; a manifold with the latter
property is {\it geometrically atoroidal}.  The product of a three-punctured
sphere with a circle is an example which is geometrically atoroidal, but not
homotopically atoroidal. The complement of an open regular neighborhood of any
torus knot is another example.

There are other results of independent interest in this paper. In particular,
Theorem~\ref{double limit theorem} is a strong general existence theorem for
limits of surface groups acting in hyperbolic three-space, and it is the main
ingredient in the proof of Theorem~\ref{mapping torus hyperbolic}. After that
proof is complete, futher results are proven concerning limits of surface
groups.  For example, Theorem~\ref{drill holes} shows how to construct
infinitely generated groups as geometric limits of surface groups of a fixed
genus. It illustrates the fact that the {\it measured} lamination is not
precisely the right structure for controlling limits of Kleinian groups.

This paper depends directly on
\refin{Theorem~5.7}{Thurston:hype1},
but it is otherwise independent of \cite{Thurston:hype1}.
Certain information about geodesic laminations, homeomorphisms
of surfaces and limits of Kleinian groups
will be assumed.  This information is summarized in
\S\ref{quasi-Fuchsian} and
\S\ref{pseudo-Anosov}.
The proofs can be gleaned from chapters 8 and 9 of \cite{GT3M}
together with \cite{FLP}, but an exposition will
also be given in parts V and VI of this series.

Theorem~\ref{mapping torus hyperbolic} was proven for the case when the fiber
is a torus minus an open disk by Troels J{\o}rgensen.

The proof of Theorem~\ref{mapping torus hyperbolic} is considerably different
from the proof of the related result for other Haken manifolds. It is
essentially a proof concerning surface groups.

The existence of a hyperbolic structure on any $3$-manifold which fibers over
the circle is paradoxical: in the universal covering space of such a
$3$-manifold, the covering space of a fiber is necessarily a uniformly bounded
distance from any image of itself by a covering transformation. This would at
first seem to be inconsistent with hyperbolic geometry: two distinct hyperbolic
planes, for instance, cannot have a uniformly bounded separation.  Two
horospheres can have uniformly bounded separation, but they have the intrinsic
geometry of the Euclidean plane, which is not possible for the universal
covering space of a fiber.

A striking image is obtained by adjoining $\sinfty$, the sphere at infinity for
hyperbolic space.  The action of the fundamental group of a hyperbolic
$3$-manifold of finite volume is {\it minimal} on $\sinfty$, that is, it admits
no closed proper invariant subsets. The closure of the universal covering of
any fiber, intersected with $\sinfty$, is a closed invariant subset, since the
fundamental group of a fiber is normal.  Consequently, the closure of the cover
of any fiber contains the entire sphere at $\infty$! (See figures 
\ref{4-1 limit}
and
\ref{figure eight limit}.)

The universal covering space of the interior of any compact surface of negative
Euler characteristic has a canonical compactification as a disk.  Here is a
more delicate fact about the fibers (\cite{Cannon:Thurston:peano}; see also
\cite{Fenley:depth1}, \cite{Minsky:peano}
and \cite{Thurston:circlesI} for generalizations
of this theorem):

\begin{figure}[htb]
\centering
\includegraphics{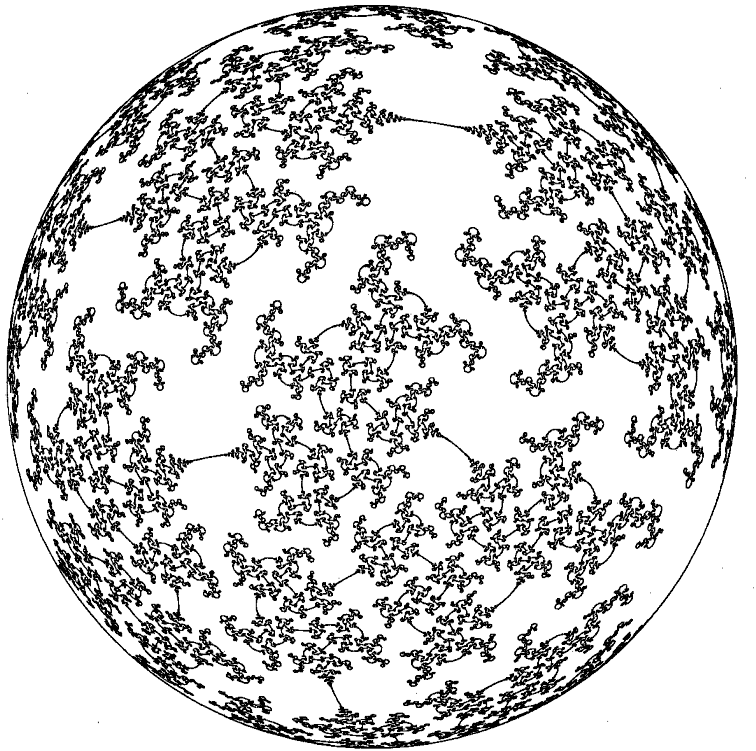}
\caption{An approximation to a sphere filling curve which arises from
the fiber of a hyperbolic three-manifold which fibers over the
circle.  The three-manifold in this case is the complement
of the figure eight knot, and the fiber is the punctured torus
bounded by the figure eight knot.}
\label{figure eight limit}
\end{figure}

\begin{theorem}[Sphere filling curve; Cannon and Thurston]\label{sphere filling curve}
Let $M^3$ be a hyperbolic $3$-manifold which fibers over $S^1$ with
fiber $F$.  Then the map $i: \tilde F \arrow \tilde M$
extends continuously to a map
$$ \bar i: D^2 \arrow D^3 .$$
The boundary of $D^2$ thus gives a sphere filling curve, or Peano
curve, on $\sinfty$.
\end{theorem}

A second purpose of this paper is to develop some of the general theory of
surface groups, especially those which share with the fibers of fibered
hyperbolic $3$-manifolds  the property that their limit set is the entire
sphere.  The result of Cannon and Thurston above is not known in the general
case.

Theorem~\ref{mapping torus hyperbolic} will finally be an easy corollary of a
general result (Theorem~\ref{double limit theorem}) which constructs limiting
surface group actions.  Precise statements along with their proofs will be
given in \S\ref{double limits}, after some background material is reviewed in
\S\ref{quasi-Fuchsian} and \S\ref{pseudo-Anosov} and a key technical theorem,
\ref{efficiency of pleated surfaces}, is proven in \S\ref{pleated surfaces}.
Theorem~\ref{mapping torus hyperbolic} will be derived in \S\ref{mapping tori}
from the result of \S\ref{quasi-Fuchsian}.

After the proof of the main theorem, there are two more sections, which take up
about half the length of this paper.  These two sections,
\S\ref{surface groups} and \S\ref{geometric limits} continue with the
analysis of the limiting behaviour of surface group actions. Some of this
information will be used later in this series, but the real thrust of this
material is toward the goal of a complete understanding of finitely generated
Kleinian groups.

The first examples of hyperbolic $3$-manifolds which fiber over a circle were
constructed and recognized by Troels J{\o}rgensen through his deep study of
limits of quasi-Fuchsian actions of the fundamental group of the punctured
torus acting on hyperbolic $3$-space.  See \cite{Jor} for a description
of some examples.  An earlier attempt I made, before I met J{\o}rgensen, to
wrestle with the question of the existence even of a Riemannian metric of
negative curvature on any such manifold, led me to really look at
homeomorphisms of surfaces and develop a geometric theory of their
classification: see \cite{Th2}, a detailed exposition of which appears
in Fathi, Laudenbach, Po\'enaru {\it et al.} \cite{FLP}.

I would like to thank Troels J{\o}rgensen for opening up this subject, and  for
sharing with me his many visions and insights. I would also like to thank
Dennis Sullivan for his relentless attack on the original unwritten,
soon-to-be-forgotten proof  of Theorem~\ref{mapping torus hyperbolic}, from
which he cut out parasitical cusps and other undesirable elements. See his
Bourbaki seminar exposition \cite{Sul2} for the first written account of
Theorem~\ref{mapping torus hyperbolic}.  Much of the further refinement of the
current version also was inspired by conversations with Sullivan.

An earlier version of this paper was circulated in preprint form in
1981.  The outline of the proof of the main theorem is substantially the
same as in the 1981 version, but the execution has been streamlined,
especially in the proof in \S3 of Theorem~\ref{pleated surfaces}.  The last
section of the current version, \S\ref{geometric limits}, is entirely new, and
its main result is new; \S6 also contains new material.
The current version is essentially the same as a version from 1986
that was refereed and accepted for publication. I dropped the ball
and did not correct and return a final copy to the journal.
This is the 1986 version, translated to \LaTeX, modified according to the
referee's comments, but with little change otherwise.


\section{Quasi-Fuchsian groups}
\label{quasi-Fuchsian}

We recall some notation from \cite{Thurston:hype1}.  Let $M^n$ be any oriented manifold and
$P^{n-1} \subset \boundary M$ be any submanifold such that the fundamental
group of each component of $P$ contains an abelian subgroup of finite index.
Then the set $\H(M, P)$ (or simply $\H(M)$ when $P = \emptyset $) is the set of
complete hyperbolic $n$-manifolds, equipped with a homotopy equivalence $f: M
\arrow N$ which sends $P$ into horoball neighborhoods of cusps of $N$.  Two
such objects are equivalent if there is an orientation-preserving isometry
between them in the homotopy class which makes the diagram commute. There are
three significant topologies on $\H(M, P)$: the {\it algebraic topology}
$\AH(M, P)$, the {\it geometric topology} $\GH(M,P)$, and the {\it
quasi-isometric topology} $\QH(M, P)$.  The maps $$ \QH(M,P) \arrow \GH(M,P)
\arrow \AH(M,P)$$ are continuous.

In the present paper we shall be dealing mainly with the case that $(M, P) =
(\bar S \cross I, \boundary \bar S \cross I)$ where $\bar S$ is a compact
surface, and $S$ will denote its interior.  An element $N \in \H(\bar S\cross
I, \boundary \bar S \cross I)$ is {\it Fuchsian} if its limit set (that is, the
limit set of the group of covering transformations of $\hy^3$ over $N$) is a
geometrical circle.  An equivalent condition is that $N$ contains a totally
geodesic oriented surface homeomorphic to $S$. The subset of Fuchsian elements
$F(S) \subset \H(\bar S\cross I, \boundary \bar S \cross I)$ is thus identified
with $\H(\bar S, \boundary \bar S)$; the three topologies restricted to this
set agree, and are the same as the Teichm\"uller space $\teich(S)$. The
classical uniformization theorem further identifies Teichm\"uller space with
the space of conformal classes of Riemannian metrics of divergence type on $S$,
that is, metrics for which every harmonic function is a constant.

An element $N \in \hs$ is {\it quasi-Fuchsian} if $N$ is quasi-isometric to a
Fuchsian manifold.  This is equivalent to the condition that the limit set for
$N$ is homeomorphic to a circle. The three topologies also agree on the set
$\QF(S)$ consisting of quasi-Fuchsian elements of $\hs$.

There is a fourth definition of a topology on $\H(M,P)$, namely the {\it
quasiconformal topology} $\QpH(M,P)$. A {\it quasiconformal map} between two
metric spaces $X$ and $Y$ is a map $f: X \arrow Y$ for which there exists a
constant $K$ such that for all $x \in X$,
$$ \lim_{r \arrow 0} \sup \left (
\frac{\sup_{x' \in S_r(x)} d(fx', fx)}
{\inf_{x' \in S_r(x)} d(fx, fx')} \right ) \le K ,$$
where $S_r(x)$ denotes the sphere of radius $r$ about $x$.

If $f$ is a quasiconformal map and if $K_1$ is a constant which works in the
above inequality for almost all $x$ on $S^2$, then $f$ is a {\it
$K_1$-quasiconformal} map.

An $\epsilon$-neighborhood of $N \in \H(M,P)$ in the quasiconformal
deformation space consists of those $N' \in \H(M,P)$ such that the two actions
of $\pi_1(M)$ on $\sinfty$ are conjugate by an orientation preserving
$\epsilon$-quasiconformal homeomorphism of $\sinfty$. The topologies $\QpH(M,P)$
and $\QH(M,P)$ are homeomorphic, by a result proved independently by Tukia
\cite{Tuk}, \cite{Tukia2}
and myself ({\it cf.} \cite{GT3M}, chapter 11 for a partial
discussion).  The current paper does not depend logically on the homeomorphism
of $\QpH(M,P)$ with $\QH(M,P)$; we will work with $\QpH(M,P)$.  Nonetheless, to
get a good intuitive feeling it is important to think about the quasi-isometric
structure of hyperbolic manifolds.

When $M$ is a complete hyperbolic manifold, let $\QpH_0(M)$ denote the component
of $M$ in $\QpH(M_0,P)$, where $M_0$ is $M$ minus horoball neighborhoods of its
cusps and $P$ is its boundary. We recall the fundamental deformation theorem
for Kleinian groups; the strong version as stated here is due to Sullivan
\cite{Sullivan1}. For any Kleinian group $\Gamma$, we denote its limit
set by $L_\Gamma$ and its domain of discontinuity by $D_\Gamma$.

\begin{theorem}[Quasiconformal deformations;
Ahlfors, Bers, \dots, Mostow, \dots, Sullivan]\label{quasiconformal deformations}
Let $M$ be any complete hyperbolic $3$-manifold with
finitely generated fundamental group.  Suppose that every component
of the domain of discontinuity of $\fund M$ is simply connected.
The conformal invariant of the quotient of the domain of discontinuity
defines a homeomorphism
$$ \mathop{conf}\colon \QpH_0(M) \arrow {D_{\fund M}}/{\fund M}.$$
\end{theorem}

Theorem~\ref{quasiconformal deformations} gives a canonical isomorphism of
$\QF(S)$ with $\teich(S) \cross \teich(S)$ since the quotient of the domain of
discontinuity for a quasi-Fuchsian group has two components, each homeomorphic
to $S$.  If $g, h \in \teich(S)$, we denote as $\qf(g, h)$ the group determined
(up to conjugacy) by $g$ and $h$ in this parametrization. Inversely, if
$\Gamma$ is a quasi-Fuchsian group, we denote the two conformal structures on
the quotient surfaces of the domain of discontinuity by $c_1(\Gamma),
c_2(\Gamma) \in \teich(S)$. The two components are distinguished by the
orientation induced from $\sinfty$: the first has positive and the second
negative orientation. (The semantic convenience of orientation as a distinction
between the two halves is one reason we are sticking for now with oriented
manifolds.)

It is easy to see that the subspace $F(S) \subset \hs$ is closed. Note that
$F(S)$ is the diagonal in the product structure for $\QF(S)$. The whole of
$\QF(S)$, on the other hand, is not closed.  For example,

\begin{theorem}[Bers slice]\label{Bers slice}
The closure in $\hs$ of any slice $x \cross F(\interior S)$
of the product structure for $\QF(S)$ is compact.
\end{theorem}

In this paper, we will find algebraic limits for sequences of quasi-Fuchsian
groups which are tending to $\infty$ not just in one factor, but in both
factors.  We will prove that such limits exists, provided the two coordinates
go to $\infty$ in directions far enough apart (thereby avoiding, for example, a
sequence of different ``markings'' of a fixed group).

The analysis will involve the geometry of the quotient manifold, so we need to
relate the conformal structure at infinity to the geometry of the
interior.  We will now give two such relations; either of these
is sufficient as a starting point for the rest of the paper.

For any element $\gamma \in \fund S$ and any $N \in \QF(S)$, let $l(\gamma)$
denote the length of the closed geodesic in $N$ homotopic to $\gamma$,
Let $\length_{+\infty}$ and $\length_{-\infty}$ denote the lengths of the closed geodesics
homotopic to $\gamma$ on the two quotient surfaces at infinity, using
their Poincar\'e metrics.

\begin{proposition}[Poincare length bounds hyperbolic length \cite{Bers2}]
\label{poincare length bounds hyperbolic length}
$$ \frac 1{l(\gamma)} \ge \frac 12 \left ( \frac 1{\length_{+\infty}(\gamma)}
+ \frac 1{\length_{-\infty}(\gamma)} \right )$$
and in particular
$$ l(\gamma) < 2 \length_{+\infty}(\gamma).$$
\end{proposition}

\begin{remark} Note that the first inequality becomes an equation when $N$ is
Fuchsian.  A sequence of examples can be constructed to show that  the
constant, $2$, in the second inequality is sharp. The idea is to ``bend'' a
Fuchsian group along a closed geodesic whose length is near zero, with the
bending angle near $\pi$. See Theorem~\ref{double limit theorem} or
\S\ref{geometric limits} for constructions which shows that there are no
inequalities of this form in the opposite sense.
\end{remark}

\begin{proof}[Proof of \ref{poincare length bounds hyperbolic length}] Let $D_+$ and $D_-$ denote the positive
and negative components of the domain of discontinuity. Then $D_+ / <\gamma>$
and $D_- / <\gamma>$ are cylinders. Calculation (say in the upper half-plane)
shows these cylinders are conformally constructed from rectangles whose
dimensions are $\pi \times \length_{+\infty}(\gamma)$ and $\pi \times
\length_{-\infty}(\gamma)$, by isometrically gluing the sides of height $\pi$.
These two annuli fit inside the torus $(S^2 - L_{<\gamma>})/<\gamma>$. This
torus is conformally constructed from a $2\pi \times l(\gamma)$-rectangle by
isometrically gluing first the two sides of length $l(\gamma)$, then the two
circles of length $2\pi$ of the resulting cylinder (with an arbitrary twist). A
standard extremal length argument (see for example \cite{AS}) gives the
proposition.
\end{proof}

A subset $A$ of a complete Riemannian manifold is {\it convex} if every
geodesic arc with endpoints in $A$ is contained in $A$.  Note that this
definition depends strongly on the topology of the ambient manifold: for
instance  a closed manifold has no proper convex subsets.  (This follows from
the fact that the geodesic flow for such a manifold is recurrent). Clearly the
intersection of any collection of
convex subsets is convex. Any non-empty convex set must
contain a possibly broken geodesic in every homotopy class. From this it
follows that a convex set in a complete hyperbolic manifold contains every
closed geodesic, since the broken geodesics in homotopy
classes $\alpha^n$, from no matter what basepoint,
wind arbitrarily near the closed geodesic in the free
homotopy class of $\alpha$, if such a closed geodesic exists.

It follows that any complete hyperbolic manifold $N$ whose
fundamental group contains at least one hyperbolic%
\footnote{We use the term 
`hyperbolic' here in its inclusive meaning, to include
the loxodromic case.}
element has a minimal non-empty
subset which is convex; this set, $C(N)$, is called the {\it convex core} of
$N$.  $C(N)$ is a codimension $0$ submanifold except in degenerate
circumstances, when it may be a submanifold of any lower dimension, possibly
with boundary.

\begin{figure}[htb]
\centering
\includegraphics[scale=.75]{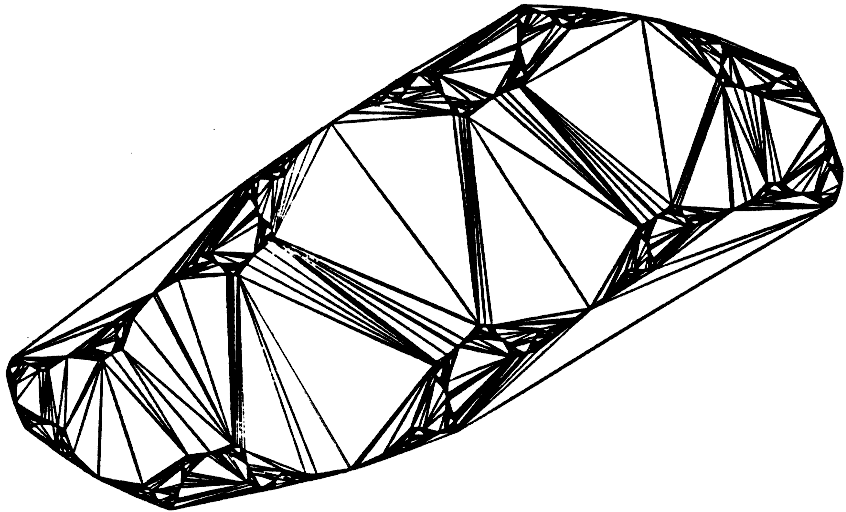}
\caption{This is an approximate drawing of the convex hull of the limit set of
a quasi-Fuchsian group. The group is a punctured torus group, that is, it is
generated by two elements whose commutator is parabolic.  This picture is in
true perspective as viewed by an observer on the sphere at infinity of
hyperbolic space. It may also be thought of as a more standard picture of the
projective ball model for hyperbolic space, where the ball is fairly large
compared to the field of vision, and the limit set has been transformed by
M\"obius transformation so that it fits within the frame.}
\label{convex hull of limit}
\end{figure}

An alternate description of $C(N)$ for a hyperbolic manifold is that it is the
quotient of the convex hull of $L_{\fund N}$
in the projective ball model for hyperbolic space
intersected with hyperbolic space and quotiented by the action of
$\fund N$.

When $N$ is $3$-dimensional and $C(N)$ is non-degenerate, the boundary of
$C(N)$ is a developable surface: it is a hyperbolic surface, homeomorphic to
the quotient of the domain of discontinuity, and isometrically embedded in $N$.
Since $\boundary C(N)$ is often not a $C^1$ submanifold,  we cannot take the
definition of an isometric embedding from standard differential geometry, and
we should clarify what definition we use: A definition which serves well is
that an embedding is isometric if every geodesic in the surface is mapped to a
rectifiable path of the same arc-length. See \cite{GT3M}, chapter 8, for
more background.

\begin{proposition}[Bounded distortion to infinity; Sullivan]
\label{bounded distortion to infinity}
There is a constant $1 < K < \infty$ such that for every complete hyperbolic
$3$-manifold $N$ and any incompressible component $S$ of $\boundary C(N)$,
there is a $K$-quasi-isometry (in the correct homotopy class) from $S$ to the
corresponding surface at $\infty$ equipped with its Poincar\'e metric.
\end{proposition}

See \cite{Sul2} for a brief proof, or \cite{EM} for a detailed
writeup which produces a concrete value for $K$.  The reasonable conjecture
seems to be that the best $K$ is $2$, but it is hard to find an angle for
proving a sharp constant.


\section{Geodesic laminations and pseudo-Anosov homeomorphisms}
\label{pseudo-Anosov}

In this section we will review some background concerning the geometry of
hyperbolic surfaces near $\infty$ in Teichm\"uller space,
and the effect of a homeomorphism of a surface on its geometry.
For details, the reader is referred to \cite{FLP}, \cite{Casson-Bleiler}
and chapters 8 and 9 of \cite{GT3M}.

We will use geodesic laminations on surfaces,
objects which are generalizations of simple closed curves
in much the same way that real numbers are generalizations of
the rational numbers.  In fact, on the torus, a simple closed curve
is described by its slope, which is a rational number, while the slope of
a geodesic lamination is a real number.

An alternative structure which serves much the same purpose
is the measured foliation.  Another very closely related object is the quadratic
differential on a Riemann surface.   
These notions are also closely related to (and partly inspired by) the work
of Nielsen on homeomoprhisms of surfaces.  See \cite{Niel-1} and
\cite{Niel-2}, and \cite{HandelThurston:Nielsen}
or \cite{Gilman} for a discussion of
Nielsen's work and its relation to mine.

Let $S$ be any complete hyperbolic surface of finite area.
A {\it geodesic lamination} $\lambda$ on $S$ is a closed subset of $S$ which is the disjoint union
of simple geodesics.  These geodesics are called the {\it leaves} of $\lambda$.
They may be either infinite simple geodesics, or simple closed geodesics.
One way to think of a geodesic lamination is to pass to $\tilde S = \hy^2$.
Each leaf of $\tilde \lambda \subset \hy^2$ is determined by its two ends
on $\cinfty$.  Therefore $\lambda$ is determined by the collection of pairs of endpoints of $\tilde \lambda$, which is a closed subset $G_\lambda$ of the open Moebius band
$M = (\cinfty \cross \cinfty - \hbox{\rm diagonal}) / \zeemod 2$,
where $\zeemod 2$ acts by interchanging coordinates.

The condition that the set $G_\lambda$ comes from a lamination is equivalent to
two closed conditions:
\begin{description}
\item[(a)] $G_\lambda$ is invariant by $\fund S$.
\item[(b)] For any two points $(a_1, a_2)$ and $(b_1, b_2)$ in $G_\lambda$, $a_1$
and $a_2$ do not separate $b_1$ and $b_2$ on $\cinfty$.
\end{description}

A geodesic lamination is defined by a subset of the projectivized tangent
bundle of a surface, namely the tangent spaces to the leaves of the
lamination.  Any Hausdorff limit of the tangent spaces of the leaves of a
family of geodesic lamination itself comes from a geodesic lamination; we call
this lamination the Hausdorff limit of the family.

Every geodesic lamination on a hyperbolic surface is a Hausdorff limit of a
sequence of geodesic laminations which have only finitely many leaves. By
forming limits of geodesic laminations with only one leaf, it is not hard to
see that if $S$ is more complicated than a $3$-punctured sphere, it has
uncountably many geodesic laminations.

A geodesic lamination on a hyperbolic surface always has zero Lebesgue measure
on $S$ (in contrast to the situation on a torus).

A {\it transverse measure} for a geodesic lamination is a measure on
$G_\lambda$ invariant by $\fund S$.  Described directly in terms of $S$, a
transverse measure assigns a measure to each curve $\alpha$ in $S$ transverse
to $\lambda$, in such a way that for any two curves $\alpha$ and $\beta$ and
any homeomorphism $f: \alpha \arrow \beta$ which takes $l \cap \alpha$ to $l
\cap \beta$ for every leaf $l$ of $\lambda$, $f$ preserves the measure. One
thinks of the measure as measuring the quantity of leaves crossed by the
$\alpha$. A geodesic lamination equipped with a transverse measure of full
support is a {\it measured lamination}.

The set of measured laminations has a topology, coming from the weak topology
on measures on the Moebius band $M$.  Denote this space, including the empty
lamination with the trivial measure, by $\ML(S)$.  The subspace of compactly
supported laminations is denoted by $\ML_0(S)$.

\begin{theorem}[ML Euclidean]\label{ML Euclidean}
$\ML_0(S)$ is homeomorphic to Euclidean space of dimension
equal to that of $\teich(S)$.
\end{theorem}

Two non-trivial measured laminations $\mu_1$ and $\mu_2$ are
{\it projectively equivalent} if their underlying laminations agree, and
one measure is a constant times the other.
The space of projective classes of measured laminations is denoted $\PL(S)$,
and the projective classes of compactly supported measured laminations is
$\PL_0(S)$.
As expected from
Theorem~\ref{ML Euclidean}, $\PL_0$ is a sphere.
Each simple closed curve defines a
point in $\PL_0(S)$.
These points are dense.

Laminations and measured laminations can easily be transferred
from one hyperbolic surface to any homeomorphic surface by means of the set $G_\lambda$,
using the fact that $\cinfty$ is a topological invariant of a surface.
Thus, the spaces $\ML(S)$, $\PL(S)$ and variations
really depend only on the
topological surface $S$.

The notion of length of a simple closed geodesic extends readily to a continuous
function
$$ \length: \teich(S) \cross \ML_0(S) \arrow \reals .$$
One way to define this extension is to define the
{\it length} $\length_S(\mu)$
of a measured lamination $\mu$ on a hyperbolic surface
$S$ to be the total mass of the ``product'' of transverse measure with
$1$-dimensional
Lebesgue measure of the leaves of $\lambda$.  More precisely, the
product measure
is defined by its rule of integration, which is an iterated integral: in
any small coordinate patch, first integrate
along the leaves of $\lambda$ with respect to Lebesgue measure,
then with respect
to the transverse measure.  When a simple closed curve is given its
tranverse counting measure, this definition agrees with the length for
the curve.

The {\it geometric intersection number} $i(\alpha, \beta) $
of two simple closed geodesics
$\alpha$ and $\beta$ is the total number of intersection points, unless $\alpha$
and $\beta$ coincide, in which case $i(\alpha, \beta) = 0$.
This also extends to a continuous function
$$ i: \ML(S) \cross \ML(S) \arrow \reals.$$
The definition for $i(\mu_1, \mu_2)$
on $S$ is defined as the total mass of a measure $\mu_1 \cross \mu_2$
on $S$, where $\mu_1 \cross \mu_2$ is the product of the two transverse
measures in any small open set where the
laminations are transverse to each other, and
zero on any leaves which the
two laminations have in common.  Transverse intersections are automatically
confined to a compact subset of $S$, so this intersection number is
always finite.  It depends only on the topological surface $S$.

\begin{theorem}[Laminations compactify Teichmuller space]\label{laminations compactify Teichmuller space}
The union $\overline {\teich(S)} = \teich(S) \cup \PL_0(S)$
has a natural topology homeomorphic to a disk.

In this topology, a sequence $\{g_i\}$ of hyperbolic
structures in $F(S)$ converges to a lamination $\mu \in \PL_0(S)$
if and only if there is a sequence
$\{\mu_i\} \arrow \infty$ of measured laminations converging projectively to
$\mu$ such that for all
$\mu' \in \ML_0(S)$ for which $i(\mu', \mu) \ne 0$,
$$ \lim_{i\arrow \infty} \frac{\length_{g_i}(\mu')}{i(\mu_i, \mu')} = 1 .$$

Furthermore, $\length_{g_0}(\mu_i) \arrow \infty$,
but $\length_{g_i}(\mu_i)$ remains bounded.

Moreover, there is a constant $C$ such that
$$ i(\mu', \mu_i) \le \length_{g_i}(\mu')
\le i(\mu', \mu_i) + C \length_{g_0}(\mu').$$
\end{theorem}

In other words, intersection number with an appropriate lamination
is quite a good approximation
to hyperbolic length near infinity in Teichm\"uller space.
It is not always the case that for a measured lamination $\overline \mu$
in the projective class of $\mu$ $\length_{g_i}(\overline \mu)$ itself
goes to zero.
For example, consider a sequence of metrics obtained by the sequence of
$i$th powers of Dehn twists about a single geodesic $\gamma$.  These metrics
converge to $\gamma$, yet $\gamma$ has constant length.  The convergence
of the $g_i$
in this case is ``tangential'' to the boundary of
Teichm\"uller space, so that a sequence of $\mu_i$ satisfying the
conditions converge
to $\gamma$ rather slowly in $PL_0(S)$.

Without the conditions in the last two paragraphs
of the theorem, the sequence of measured laminations
$i \length_{g_0}(\gamma) \gamma$ would serve. 

If $\beta$ is any simple closed curve which intersects $\gamma$, then a certain
constant multiple of the sequence of images of $\beta$ by the $i$th power
of the inverse Dehn twist will satisfy the conditions for $\mu_i$ except for
the last paragraph. 

The laminations $\mu_i$ are constructed in \cite{Th5}, in
the course of development of cataclysm coordinates. This is quite related to
the proof of Theorem \ref{efficiency of pleated surfaces} as well; in cataclysm
coordinates, given a lamination $\lambda$ and a metric $g$, a measured
foliation $\mathcal F$ is constructed such that the length along $\lambda$ in $g$
exactly agrees with its transverse measure.  The first stage of the polygonal
approximations made in the proof of \ref{efficiency of pleated surfaces} show
that the measured lamination $\nu$ defined by $\mathcal F$ serves well to estimate
length.

\medskip

Note that if the $\mu_i$ are normalized to converge to $\mu$ as
measured laminations, rather than simply projective
laminations, then $\length_{g_i}(\mu_i) \arrow 0$.

\bigskip

Pairs of laminations are often important in the theory of surfaces.  A pair
$(\mu, \nu)$ of laminations is called {\it binding} if
\begin{description}
\item[(a)] they have no leaves in common, and
\item[(b)] for each component $U$ of the complement of the union of $\mu$ and
$\nu$, the metric completion of $U$ is
either a compact polygon, 
a polygon with one ideal vertex which tends to a cusp, or
a punctured polygon, where the puncture is a cusp of the surface.
\end{description}

\begin{proposition}[Lamination crosses binding]\label{lamination crosses binding}
A pair of laminations $(\mu, \nu)$ on a surface $S$ is binding if and only if
every simple geodesic on $S$ has at least one transverse
intersection with a leaf of $\mu$ or a leaf of $\nu$.

A pair of compactly supported measured laminations $(\mu, \nu)$
is binding if and only if for every compactly supported measured
lamination $\lambda$,
$$ i(\lambda, \mu) + i(\lambda, \nu) > 0$$
\end{proposition}

\begin{proof}[Proof of \ref{lamination crosses binding}]
If $(\mu, \nu)$ is binding, then every leaf of $\mu$ crosses at least one
leaf of $\nu$ (and {\it vice versa}).
For suppose, on the contrary, that $(\mu, \nu)$ is a pair of laminations
such that $\mu$ has a leaf $l$ which doesn't meet $\nu$.
Then the closure of $l$ is a sublamination $\lambda$
of $\mu$; if $\lambda$ meets $\nu$, then the two laminations have leaves
in common, so the pair is not binding.
Otherwise, since $\nu$ is closed, $\lambda$ has a neighborhood
not meeting $\nu$.
This  violates (b), so the pair is not binding.

If $(\mu, \nu)$ is binding, and if $g$ is any simple geodesic
on the surface, then if $g$ intersects some component of the
complement of the union of $\mu$ and $\nu$, it is clear that
$g$ intersects one of the two laminations transversely, by (b).
The only other possibility is that $g$ is a leaf of $\mu$ or of $\nu$,
in which case it intersects the other lamination transversely,
by the preceding paragraph.

On the other hand, suppose that $(\mu, \nu)$ is a pair of
laminations such that every simple geodesic intersects one or
the other transversely.  Then, in particular, they have no leaves in common.
If any region of the complement is not simply-connected, then
its fundamental group must be parabolic, otherwise it would contain
a simple closed geodesic.  The metric completion
of a region of the complement has boundary
made up of segments of leaves of $\mu$ and of $\nu$, so it is polygonal
with possibly some missing vertices.  There can be no more than
one missing vertex, for otherwise, there would be a simple geodesic
connecting one ideal vertex to the other. The neighborhood of
any ideal vertex
must tend toward a cusp of $S$, because if it recurred in
compact subsets of $S$, one could form a limit, and construct
a new simple geodesic with no transverse intersections with $\mu$
or $\nu$.
If the boundary
has any missing vertices, then the region is simply-connected, or again
there would be a simple geodesic tending to the ideal vertex at either end,
looping around the fundamental group. 

The only remaining possibilities are
those mentioned in (b).
This establishes the claim of the first paragraph of the proposition.

Now let us suppose that $\mu$ and $\nu$ are compactly supported measured
laminations.  From the previous condition, it
follows that if $\lambda$ is a non-trivial measured
lamination such that $i(\lambda, \mu) = i(\lambda, \nu) = 0$,
then $(\mu, \nu)$ is not binding.

To establish the converse, recall (for instance, from
\refin{Proposition 5.4}{Thurston:hype1}) that
every measured lamination is the finite union of minimal sublaminations.
Therefore, if $\mu$ and $\nu$ have any leaves in common, the set of
common leaves constitutes a finite union of minimal sets of both.
It inherits a transverse measure say from $\mu$, to give it the
structure of a measured lamination $\lambda$ for which
$i(\lambda, \mu) = i(\lambda, \nu) = 0$.

If a region of the complement of $\mu \cup \nu$ violates condition (b),
there are two possibilities.  If the region has a fundamental
group not consisting of parabolic elements, there is a measured lamination
$\lambda$ supported on a simple closed geodesic and not intersecting
$\mu$ or $\nu$.  If the boundary of the region has an ideal vertex,
it cannot tend to a cusp, since $\mu$ and $\nu$ have compact support.
Therefore, its limit set gives a minimal set in $\mu$ or $\nu$,
or both.  It inherits a transverse measure from $\mu$ or from $\nu$, 
thereby defining a measured lamination having intersection number
0 with $\mu$ and $\nu$.
\end{proof}

For later use, we need:

\begin{proposition}[Binding confinement]\label{binding confinement}
Let $\mu_1$ and $\mu_2$ be any two measured laminations which have no leaves in
common, and for which
$S \setminus (\mu_1 \cup \mu_2)$ consists of pieces which are
simply connected or neighborhoods of cusps.

 Then for any constant $C > 0$, the set of hyperbolic structures $g$ on $S$
for which 
$$ \length_g(\mu_1) \le C$$
and
$$ \length_g(\mu-2) \le C$$
is compact.
\end{proposition}

\begin{proof}[Proof of \ref{binding confinement}]
Two such laminations
$\mu_1$ and $\mu_2$ have the property that for any measured lamination
$\nu$, $i(\mu_1, \nu) > 0$ or $i(\mu_2, \nu) > 0$.  Applying this to the
case that $\nu$ ranges over a sequence $\mu_i$ from
\ref{laminations compactify Teichmuller space},
we see that no sequence of metrics for which both $\mu_1$ and $\mu_2$
have bounded lengths can approach the boundary of Teichm\"uller space.
\end{proof}

For a different, directly geometric proof of a more general proposition, see
Theorem \ref{laminations bind subsurface}. Here is another fact for future
reference:
  
A homeomorphism $\phi: S \arrow S$ is
{\it reducible} if $\phi$ permutes some finite system of
disjoint simple closed geodesics.  The terminology is
based on the fact that by cutting $S$ along such a system
of curves, one can reduce the study of $\phi$ to the study of
a homeomorphism of a simpler surface.  This is analogous
to the notion of a reducible three-manifold, or (closer to home)
a torus-reducible three-manifold.

\begin{theorem}[Classification of surface homeomorphisms]\label{classification of surface homeomorphisms}
Every homeomorphism $\phi'$ of $S$ is isotopic to a homeomorphism
$\phi$ which either
\begin{description}
\item[(i)] has finite order
\item[(ii)] is reducible, or
\item[(iii)] does not satisfy (i) or (ii), and
preserves a unique pair of projective classes of measured
laminations, $\mu_s$ and $\mu_u$.
\end{description}

In the case (iii), for any point $x \in \overline{\teich(S)}$, if $x \ne \mu_s$
and $x \ne \mu_u$, then
$$ \lim_{n\arrow +\infty} \phi^n(x) = \mu_s$$
and
$$ \lim_{n\arrow -\infty} \phi^n(x) = \mu_u.$$
There is a constant $\lambda > 1$ such that $\phi$ multiplies the
transverse measure for $\mu_u$ by $\lambda$ and the transverse
measure for $\mu_s$ by $1/\lambda$.
\end{theorem}

Case (iii) is both most common and most interesting.
The isotopy class of $\phi$ in case (iii) is called
{\it pseudo-Anosov}.  There is a more geometric form for $\phi$,
in this case, if the measured laminations are replaced by
measured foliations.  When this is done, $\phi$ becomes
almost an Anosov homeomorphism: it is Anosov in the complement
of a finite collection of points.

From any homeomorphism $\phi$ of a surface $S$ to itself, a $3$-manifold
called the {\it mapping torus} of $\phi$ can be constructed
by first forming the product $S \cross I$, and then gluing
$S \cross \{1\}$ to $S \cross \{0\}$ using $\phi$.  The resulting
3-manifold $M_\phi$ fibers over $S^1$,
and $\phi$ is called the monodromy of the fibration of $M_\phi$.

Actually, only the isotopy class of $\phi$ is determined
by the fibering $M_\phi \arrow S^1$, and the isotopy
class of $\phi$ determines $M_\phi$ up to
homeomorphisms which commute with $M_\phi \arrow S^1$.

\begin{proposition}[Non-pseudo-Anosov mapping torus]\label{non-pseudo-Anosov mapping torus}
If $\phi$ is a homeomorphism of a compact surface $\bar S$
of negative Euler
characteristic, then
\begin{description}
\item[(i)] $\phi$ is isotopic to a homeomorphism of finite order if
and only if the interior of $M_\phi$
admits a complete $\hy^2 \cross \reals$ structure of finite volume
\item[(ii)] $\phi$ is isotopic to a reducible homeomorphism if and only if
$M_\phi$ admits an embedded incompressible torus not isotopic to $\boundary M$,
and
\item[(iii)] otherwise $\phi$ is pseudo-Anosov.
\end{description}
\end{proposition}

\begin{proof}[Proof of \ref{non-pseudo-Anosov mapping torus}]
If $\phi$ is isotopic to a homeomorphism of finite order, then there
is a hyperbolic structure on $S$ which is invariant by $\phi$ up to isotopy.
From this, a $\hy^2 \cross \reals$ structure is constructed on $M_\phi$.

If the interior of
$M_\phi$ admits a $\hy^2 \cross \reals$ structure, that structure
defines on $M_\phi$ a codimension two transversely hyperbolic foliation.
It was shown in \cite{GT3M} that such a foliation either describes
a Seifert fibration, or the manifold fibers over the circle with fiber
a torus.  Clearly the former case obtains.  Consider the projection
of the fiber $\bar S$ of the fibration over $S^1$ to the base of the
Seifert fibration.  There is an induced map from the fundamental
group of $\bar S$ to the fundamental group of the base orbifold.  The image
is a normal subgroup.  The quotient group must be finite, since there
is no finitely generated normal subgroup of the fundamental
group of a hyperbolic $2$-orbifold with quotient group $\integers$.
The map is also injective on the level of fundamental group.
This implies it is homotopic to a covering map over the
base orbifold.  Then $\phi$ is isotopic to a deck transformation,
so it has finite order.

If $\phi$ is isotopic to a reducible homeomorphism, then the trace of
a reducing curve, dragged around by $\phi$ until it comes back to itself,
with the same orientation,
gives an incompressible torus.

Conversely, if there is an incompressible torus, consider the
projection of the torus to the circle.  It can be perturbed
to have Morse-type singularities.
A standard $3$-manifold argument shows that
the singularities can all be eliminated by an isotopy of the torus,
so that the torus becomes transverse to each fiber.
The intersection with one of the fibers defines a family of
reducing circles for some homeomorphism isotopic to $\phi$.
\end{proof}


\section{An estimate for the shapes of certain pleated surfaces}
\label{pleated surfaces}

Let $\bar S$
be a compact surface and $N \in \hs$ be a hyperbolic $3$-manifold.
For purposes of reference, fix a complete hyperbolic structure of finite
area on $S$.  Most non-trivial simple closed curves (except
parabolics) are realized as closed geodesics in $N$.  Similarly,
most geodesic laminations on $S$ have realizations in $N$.
In particular, the realizable laminations contain an open dense
subset of $\PL(S)$.
This theory is developed in chapters 8 and 9 of \cite{GT3M}.
In the special case $N \in \QF(S)$, all geodesic laminations on
$S$ are realized in $N$.

The quantity $\length_N(\mu)$ of a realizable measured lamination
$\mu$ in a hyperbolic $3$-manifold $N$ is defined similarly
to the length of a lamination on a surface, as the total
mass of a ``product'' measure formed from
1-dimensional Lebesgue measure
along the leaves with transverse invariant measure. To put it
another way, in a local coordinate system, the transverse invariant
measure is a measure on the space of leaves of $\mu$; the integral
of a continuous function $f$ with respect to the product
measure is defined by first integrating $f$ with respect to
Lebesgue measure on each local leaf, then integrating with respect
to the transverse measure.  A further complication is that
$\mu$ is only mapped into $N$, but not necessarily embedded:
one should do the computation in the domain, by
pulling back Lebesgue measure from $N$.

\begin{proposition}[Length continuous]\label{length continuous}
The length of laminations is a continuous function on the set
$R \subset \PL(S) \cross \hs$ of realizable laminations.
\end{proposition}

The details of this and other basic facts about measured laminations in two and
three dimensions will be proven in part VI of the series.

A geodesic lamination $\lambda$ of a surface $S$ is {\it maximal} if it is
maximal among all geodesic laminations, which is equivalent to the condition
that each region of $S - \lambda$ is isometric to the interior of an ideal
triangle.  Every lamination is contained in a maximal lamination, although a
lamination may be maximal among measured laminations without being maximal in
the current sense: to extend it to a maximal lamination, one must often make do
with a lamination admitting no transverse measure of full support. Every
surface admits maximal laminations which have only finitely many leaves.  See
for example \refin{Figure 2.1}{Thurston:hype1} for a way to construct laminations by
``spinning'' triangulations.

If $\boundary \bar S \ne \emptyset$, the situation can be simplified by
allowing no closed leaves on $S$. Then $\lambda$ becomes the $1$-skeleton of an
``ideal triangulation'' of $S$.

A lamination $\lambda$ with a finite number of leaves is realizable provided
all its closed leaves are realizable.  If $\lambda$ is maximal, then it
determines a pleated surface
$$ f_\lambda: P_\lambda \arrow N ,$$
where $P_\lambda$ is a complete hyperbolic structure on $S$
and $f_\lambda$ is an isometric map of $P_\lambda$ into $N$ which
folds or pleats, at most, along $\lambda$.

If $\lambda$ is a maximal lamination with only finitely many leaves and $\mu$
is a measured lamination of compact support, we will define a quantity
$$ a(\lambda, \mu)$$
which in some sense measures the complexity of $\mu$
relative to $\lambda$.  We will define $a(\lambda, \mu)$ first in the case
that $\mu$ is a simple closed curve.  If $\mu$ is a leaf of
$\lambda$, then we define $a(\lambda, \mu) = 0$.
Otherwise, there is at most a countable set of intersections of the
two laminations.  If $x$ and $y$ are two intersection points with
no intervening intersections along $\mu$, then the leaves of $\lambda$
through $x$ and $y$ are two sides of an ideal triangle of $\lambda$ (since
$\lambda$ is maximal, by hypothesis).  These leaves are therefore
asymptotic on either one side of $\mu$, or the other.
If $x$, $y$, and $z$ are three successive intersection points along $\mu$,
and if the leaves of $\lambda$ through $x$ and $y$ are asymptotic on the
opposite side of $\mu$ to the leaves of through $y$ and $z$, then let us
call $y$ a {\it boundary intersection}.

\begin{figure}[htb]
\centering
\includegraphics{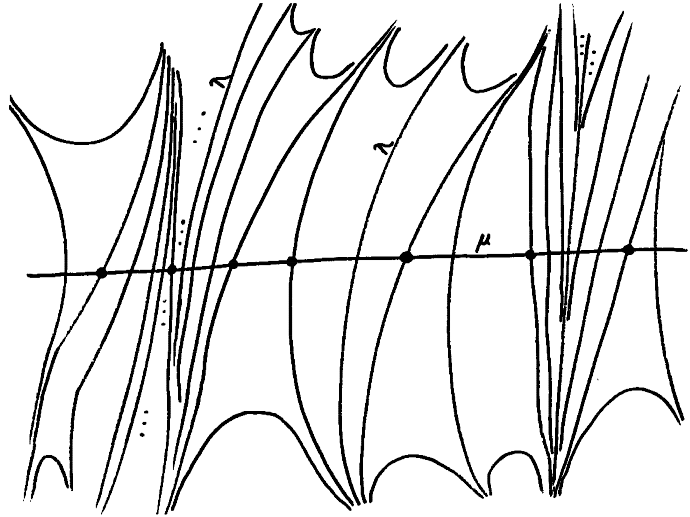}
\caption{The boundary intersections of a geodesic $\mu$ with a lamination
$\lambda$ are marked in the universal cover, $\hy^2$.  Even where leaves of
$\lambda$ accumulate, the boundary intersections are isolated.}
\label{boundary intersections and alternation number}
\end{figure}

It can also happen that a point $y$ of intersection of $\lambda$ and $\mu$
is an accumulation point of leaves of $\lambda$.  In such a case,
the leaf of $\lambda$ through $y$ is a closed leaf, and the leaves of
$\lambda$ spiral to the closed leaf
on both sides.  We will call such an intersection point
a {\it boundary intersection}
if the directions of spiraling on the two sides are different.
Note that both kinds of boundary intersection are isolated, and so
there is a finite number of boundary intersections in all.
We define
$a(\lambda, \mu)$ (in the case that $\mu$ is a simple closed curve)
to be the total number of boundary intersections of $\lambda$ and $\mu$.
In other words, $a(\lambda, \mu)$ measures the total number of times the
direction of asymptoticity of the leaves of $\lambda$ changes as you go
around $\mu$.

When $\mu$ is a general measured lamination, the situation is much the
same. Designate a point $y$ of transverse intersection of $\lambda$ and
$\mu$ a {\it boundary intersection} if the direction of asymptoticity
of the leaves of $\lambda$ changes at $y$.  Then compute $a(\lambda, \mu)$
as the total $\mu$-transverse measure of the set of boundary
intersection points.

\begin{proposition}[Continuity of alternation number]\label{continuity of alternation number}
For a fixed maximal finite lamination $\lambda$, the alternation number
$a(\lambda, \mu)$ is a finite-valued,
continuous function of the measured lamination $\mu$.
\end{proposition}

\begin{proof}[Proof of \ref{continuity of alternation number}]

Associated with a measured lamination $\mu$ is a measure $M(\mu)$ on the
projective bundle $\proj(S)$ of $S$, defined as the product
of 1-dimensional Lebesgue measure with the transverse measure.
We will construct a continuous function $C(\lambda)$ such that for all
$\mu$,
$$ a(\lambda, \mu) = \int_{\proj(S)} C(\lambda) d M(\mu) .$$
This implies that $a(\lambda, \mu)$ is continous in $\mu$.
 
The idea for constructing $C(\lambda)$ is that instead of counting
boundary intersections directly, we can spread the contribution out
so that it becomes an integral
over a larger portion of $\proj(S)$.

Every isolated leaf $l$ of $\lambda$, lifted to the
universal cover of the surface, separates two ideal triangles.
The union of the two triangles forms an ideal quadrilateral.  If
a geodesic $m$ crosses $l$, it is a boundary intersection iff $m$ passes
through opposite sides of this quadrilateral.

Choose a continous
function on the unit interval with integral $1$ which
is 0 at the two endpoints of the interval.  By scaling,
this transfers to a function whose integral is $1$ on
an arbitrary interval.

For each quadrilateral obtained by removing an isolated leaf of $\lambda$ in the
universal cover, and for each geodesic $m$ which intersects opposite sides of
the quadrilateral, scale this function to the intersection
of $m$ with the quadrilateral, and lift to
$\proj(\hy^2)$.  Define a function in $\proj(\hy^2)$ by adding up over
all quadrilaterals and all geodesics.
For any point of $\proj(\hy^2)$, there are at most $3$
contributions, so this gives a well-defined function on $\proj(\hy^2)$.  It is
continuous, because when a sequence of geodesics $m$ which cross opposite
sides of a quadrilateral converge to one that doesn't, the length goes to
$\infty$ so the function tends to $0$.
This function is invariant by deck transformations, so it gives a continous
function on $S$.

Similarly, if $\alpha$ is a closed geodesic
of $\lambda$ for which the spiraling is in opposite directions on its two sides,
we add a
contribution for each geodesic $m$ which crosses $\alpha$,
supported on the intersection of $m$ with
an $\epsilon$ neighborhood of $\alpha$.  The result, after adding all the
contributions, is a continuous function $C(\lambda)$, with the desired
properties.

\end{proof}

The following theorem gives an estimate of how efficiently certain homotopy
classes can be represented on a pleated surface with finite pleating
locus:

\begin{theorem}[Efficiency of pleated surfaces]
\label{efficiency of pleated surfaces}

Let $\bar S$ be a fixed compact surface, possibly with boundary.
For any $\epsilon > 0$ there is
a constant $C < \infty$ such that the following holds:

Let $\lambda$ be any finite maximal lamination on $S$.

Let $N$ be any element of $\hs$.
Suppose that no closed leaf of $\lambda$
has length less than $\epsilon$ in $N$, so that in particular,
no such leaf is parabolic, and consequently a surface
$f_\lambda: P_\lambda \arrow N$ exists which is pleated along $\lambda$.

Let $\mu \in \ML_0(S)$ be any compactly-supported measured lamination
which is realizable in $N$.

Then
$$ \length_N(\mu) \le \length_{P_\lambda}(\mu)
\le \length_N(\mu) + C a(\lambda, \mu) .$$
\end{theorem}

This inequality does not generalize to an inequality for lengths
on more general pleated surfaces.  The finite combinatorial complexity
of $\lambda$ is crucial.

When $\mu$ is a not necessarily
realizable lamination in $\ML_0(S)$,
define $\llength_N (\mu)$ to be the lim inf of lengths of nearby
realizable laminations.  If $\mu$ has only one component, then it follows
from \cite{GT3M} or \cite{Bon} that $\llength_N(\mu) = 0$.

\begin{corollary}[Efficiency for unrealizable laminations]\label{efficiency for unrealizable laminations}
Theorem \ref{efficiency of pleated surfaces} works for arbitrary
$\mu \in \ML_0(S)$ if
$\length_N(\mu)$ is replaced by $\llength_N(\mu)$.
\end{corollary}

\begin{proof}[Proof of \ref{efficiency for unrealizable laminations}]
Apply the theorem to an appropriate sequence of realizable
laminations converging to $\mu$.
\end{proof}

\begin{proof}[Proof of Theorem \ref{efficiency of pleated surfaces}]

In view of proposition \ref{length continuous},
it suffices to prove \ref{efficiency of pleated surfaces}
in the case that $\mu$ is a simple closed curve, since simple
closed curves are dense among realizable laminations
and the inequality is homogeneous in $\mu$.

The idea of the proof is to represent a simple closed curve
$\mu$ by a polygonal path on $P_\lambda$ whose number of
sides is $O(a(\lambda, \mu))$ \def\alm{a(\lambda, \mu)} and
which follows along the leaves of $\lambda$ except for $O(\alm)$
of its
length.  By an application of the uniform injectivity theorem
\cite[Thm. 5.7]{Thurston:hype1}, one finds that
this polygonal path on $P_\lambda$, when mapped to $N$,
cannot ``double back'' very much, so that it cannot be too
inefficient.

From now on in the proof, $\mu$ will be a simple closed geodesic on
$P_\lambda$.
We shall homotope $\mu$ to a polygonal curve.

We may assume that $\mu$ is not
a leaf of $\lambda$, for in that case there is nothing to prove: its length
in $N$ equals its length in the pleated surface.

We may also assume that $\mu$ is not a geodesic shorter than
$\epsilon$ on $S$, for if it is not a leaf of $\lambda$, then
$a(\lambda, \mu) \ge 2$, and the inequality is trivial.

Otherwise, $a(\lambda, \mu) > 0$, and there is a chain of length
$a(\lambda, \mu)$ consisting of the leaves
of $\lambda$ on which there are boundary intersections
with $\mu$.
Successive leaves in the chain are asymptotic, and if the asymptotic
pairs of ends of leaves are replaced by short jumps, one obtains
a curve homotopic to $\mu$.

It is easiest at this point to describe the picture in the
universal cover of the surface, $\hy^2$.  Let $R$ be a strip
of constant width, say $.5$, about the geodesic $\tilde \mu$.
For concreteness, we can consider the picture in the Poincar\'e disk
model, with $\tilde \mu$ being the horizontal diameter of the
disk.

We have an infinite chain
$g_i \  [i=-\infty \dots \infty]$ of geodesics connecting the
two endpoints of $\tilde \mu$, with endpoints of
successive geodesics meeting at $\cinfty$.  Each $g_i$ crosses
$R$.  Let $a_i$ be the lower
endpoint of $R \cap g_i$, and $b_i$ the upper endpoint.  Let $A_i$ and
$B_i$ be the
(hyperbolic) perpendicular projections of $a_i$ and $b_i$ to $\tilde \mu$.
Then $A_i < A_{i+1}$, $B_i < B_{i+1}$, and $|A_{i+1} - B_i|$ is
bounded by an {\it a priori} constant.

If $B_i > A_{i+1}$, let
$x_{i+1}$ be the point of $g_{i+1}$ whose perpendicular projection
to $\tilde \mu$ is $B_i$; otherwise, let $x_{i+1} = a_{i+1}$.

\begin{figure}[htb]
\centering
\includegraphics{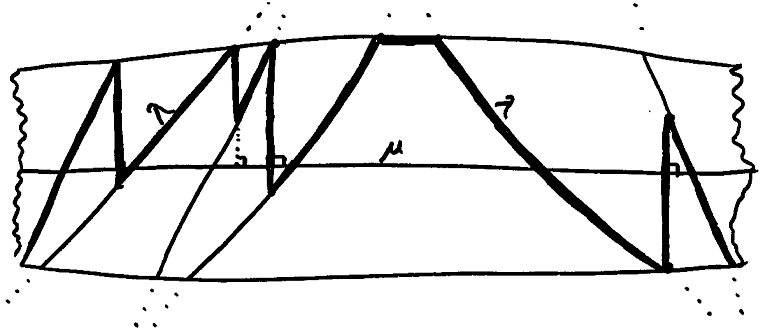}
\caption{An arbitrary simple closed geodesic $\mu$ on a pleated surface can be
approximated by a polygonal path with $2a(\lambda,\mu)$ sides, which follows
leaves of $\lambda$ except on $O(a(\lambda, \mu))$ of its length, and whose
length exceeds that of $\mu$ only by $O(a(\lambda, \mu))$. The path simply
follows appropriate segments of the sequence of asymptotic geodesics given by
its boundary intersections with $\lambda$.  The possibly long segments of this
polygonal path which follow leaves of $\lambda$ are mapped efficiently in any
$\lambda$-pleated surface.  By making further adjustments to avoid problems
with the thin parts of the surface, the polygonal
path can be forced to map efficiently as a whole on a $\lambda$-pleated
surface, up to an additive constant which is a bounded multiple of $a(\lambda,
\mu)$.}
\label{polygonal approximation for pleated surface}
\end{figure}

Define a polygonal path $\tilde m \subset \hy^2$ to consist of
the the geodesic segments
$$ \dots \  \overline{x_{i} b_{i}} \  \overline {b_{i} x_{i+1}} \  \dots .$$
This polygonal path consists of intervals of the leaves $g_i$,
interspersed with jumps of bounded length from one leaf to another.
The projection to $P_\lambda$ is a polygonal curve $m$ whose length exceeds
the length of $\mu$ only by a constant times $a(\lambda, \mu)$.
Furthermore, $m$ lies on leaves of $\lambda$ except for a
portion of its length which is a constant times $a(\lambda, \mu)$.

Next, we will make a homotopy of $m$ to a new curve $n$ with similar
properties to $m$, but which also has the property that none of the
jumps between leaves of $\lambda$ occur in the thin part of $P_\lambda$.
We may assume that $\mu$ is not a short geodesic
The curve $m$ enters only those thin parts of $P_\lambda$ which are neighborhoods
of short geodesics.  Since, by hypothesis, no closed
geodesics of $\lambda$ are short, each leaf of $\lambda$
which enters such a thin part exits on the opposite side.
All leaves of $\lambda$ which enter the thin part remain
fairly close together and roughly parallel for the entire intersection
with the thin part.

For any interval of $m$ which crosses a thin part of $P_\lambda$, 
we can inductively replace the first segment on a leaf in the
thin part and the first jump in the thin
part by a jump just outside the thin part, and a longer segment on
a leaf.  Each such step increases the length of the curve by at most
a constant. The total number of steps is at most $a(\lambda, \mu)$, so
in the end we obtain a polygonal curve $n$ with the desired properties.

\medskip
The image $P_\lambda (n)$ of $n$ in the three-manifold $N$ is probably not
polygonal, but we can homotope it to a polygonal path $p$ by homotoping the
image of each jump segment of $n$ to its geodesic in $N$.  The length of
$p$ differs from the length of $\mu$ by at most
a constant times $a(\lambda, \mu)$. The polygon $p$ has the same number
of sides as $n$: $2 a(\lambda, \mu)$.

Let $\mu'$ be the geodesic of $N$ homotopic to $P_\lambda(\mu)$.
Construct a pleated annulus $A$ representing the homotopy from $p$
to $\mu'$ by spinning around $\mu'$:
concretely, for each side of $p$, form a triangle based on that side
with third ``vertex'' spiralling infinitely around $\mu'$.
Form $A$ by gluing all these triangles together, and completing
with $\mu'$.
The area of each triangle
is less than $\pi$, so the area of $A$ is less than $2\pi a(\lambda,\mu)$.

In the intrinsic hyperbolic metric of $A$, the portion of the $p$ boundary
component of $A$ which admits a regular neighborhood of width
$\epsilon$ in $A$ is bounded by $area(A)/\epsilon$, which is
less than a constant times $a(\lambda, \mu)$.
Any remaining length of $p$ consists of segments on leaves of
$\lambda$ in $N$ which run close
and nearly parallel to some other portion of the boundary of $A$.

We claim that the nearby portion of the boundary of $A$ is part of
$\mu'$, not part of $p$, except possibly for a total length bounded by
a constant times $\alm$.

\begin{figure}[htb]
\centering
\includegraphics[scale=.75]{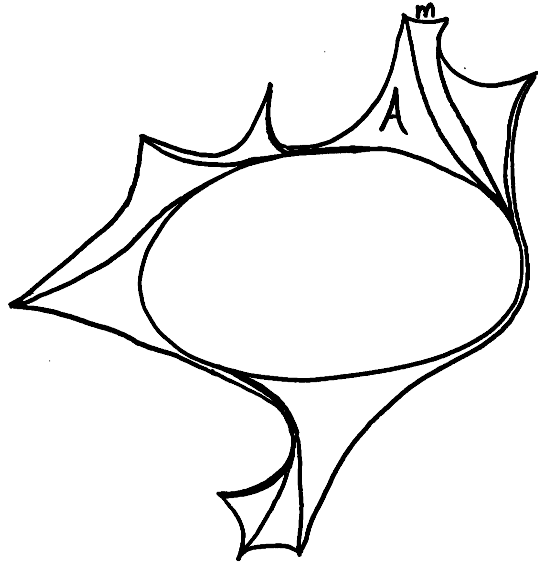}
\caption{A pleated annulus $A$ constructed between a polygonal curve $m$ and
its closed geodesic has area bounded by the number of sides.  If the polygonal
boundary component is large compared to the number of sides, this means that
much of its length must be nearly parallel to another portion of the boundary
of $A$.  We have constructed the polygon so that the length of $m$ is in excess
of a constant times $a(\lambda, m)$ consists of leaves of $\lambda$.  The
leaves of $\lambda$ are forced apart in light of the  the uniform injectivity
theorem so the lengths of the two boundary components differ by only
$O(a(\lambda, m))$.}
\label{pleated annulus}
\end{figure}

For any portion of $p$ which is not in $P_{\lambda,\thin}$,
this will follow fairly
directly from \refin{Theorem 5.7}{Thurston:hype1}. In fact, if there is any short arc
$\alpha$ on $A$ connecting two segments of $p$ on leaves of $\lambda$ but not
in $P_{\lambda,\thin}$, $\alpha$ would be homotopic to a short arc on $S$, by the uniform
injectivity theorem.  From the picture on $A$, $\alpha$ would also be homotopic
to an interval of $p$ and of $n$, so it gives a way to shorten $n$. This is
possible at most for a part of $n$ of length less than a constant times
$a(\lambda, \mu)$.

The polygonal curve $n$ on $P_\lambda$ has the property that if there is any interval
$I$ of $n$ with endpoints in $P_{\lambda,\thin}$ which is homotopic rel endpoints to
$P_{\lambda,\thin}$, then this interval is contained on a single leaf of $\lambda$.  The
map from components of the thin set of $P_\lambda$ to the thin set of $N$ is
injective.  Therefore, two distinct segments of $p$ which are in the thin part
of $P_\lambda$ cannot be close together on $A$.

Therefore, all but a constant times $a(\lambda, \mu)$ of the length of $p$ runs
close to $\mu'$ on $A$. This implies that
$$ \length_{P_\lambda}(\mu)  \le \length(p) \le
\length_N(\mu) + C a(\lambda, \mu),$$
as desired.
\end{proof}

Theorem \ref{efficiency of pleated surfaces} is false without the stipulation
that closed leaves of $\lambda$ have length at least $\epsilon$ in $P_\lambda$,
or equivalently, in $N$.  It is possible to construct examples where there is a
pleated surfaces which has a nearly $180 \degrees$ fold along a very short
closed leaf of $\lambda$.  Any geodesic of $P_\lambda$ which intersects such a
fold can be shortened considerably in $N$.

The limiting case of complete inefficiency is for surface groups which have an
accidental parabolic. Take $\lambda$ to be a lamination having a closed leaf
whose homotopy class is parabolic.  Then there is no actual pleated surface
$P_\lambda$, but there is still a type of ``pleated surface with nodes'', which
goes off to $\infty$ at the parabolic curve.  Any homotopy class passing
through the parabolic curve has infinite length on $P_\lambda$, but finite
length in $N$.

Here is a slightly stronger version of the theorem, which might be useful
sometime:

\begin{theorem}[Restricted efficiency of pleated surfaces]\label{restricted efficiency of pleated surfaces}

Let $\bar S$ be a fixed compact surface, possibly with boundary. For any
$\epsilon > 0$ there is a constant $C < \infty$ such that the following holds:

Let $\lambda$ be any finite maximal lamination on $S$.

Let $N$ be any element of $\hs$. Let $\lambda_1$ be the measured lamination
which is the sum of the closed leaves of $\lambda$ of length less than
$\epsilon$ in $N$, and let $\lambda_2$ be the union of closed leaves of length
zero in $N$, that is, closed leaves which are parabolics.

Then there is a pleated surface $f_\lambda: P_\lambda \arrow N$
pleated along $\lambda \setminus \lambda_2$,
where $P_\lambda$ is a complete hyperbolic
structure on $S - \lambda_2$.  
For any $\mu \in \ML_0(S)$ such that $i(\mu, \lambda_1) = 0$,
$$ \length_N(\mu) \le \length_{P_\lambda}(\mu)
\le \length_N(\mu) + C a(\lambda, \mu) .$$
\end{theorem}
\begin{proof}[Proof of \ref{restricted efficiency of pleated surfaces}]
The proof of \ref{efficiency of pleated surfaces} works also for this
statement. The absence of short geodesics of $\lambda$ was used only in the
construction of a polygonal approximation to the curve $\mu$, such that no
jumps take place inside $N_\thin$.  The difficulty is averted if $\mu$ does not
enter the components of $P_{\lambda,\thin}$ surrounding short closed leaves of
$\lambda$.  But any {\it simple} closed geodesic remains entirely outside such
a thin neighborhood, unless it crosses it.
\end{proof}

\begin{remark} In fact, both
Theorem \ref{efficiency of pleated surfaces} and the more
general version, Theorem \ref{restricted efficiency of pleated surfaces}
can be extended to curves on the surface which are not simple, and to
more general measured laminations  tangent to the
geodesic flow in $\proj(S)$, which do not
project to simple laminations on $S$.  The alternation number makes
just as much sense, and the polygonal approximations work equally well.
A non-simple geodesic may enter the neighborhood of a short geodesic
or a neighborhood of a cusp, wind around many times, and then exit
through the boundary component of the thin set where it entered.
In this case, one can again push all the jumps between leaves of
$\lambda$ in a polygonal approximation out of the thin set, one
by one, adding only a bounded multiple of the alternation number.
\end{remark}


\section{Existence of double limits of quasi-Fuchsian groups}
\label{double limits}

Using the results of
\S\ref{pleated surfaces},
we can now prove a good existence theorem for limits of surface group
actions in $\hy^3$  In particular, the theorem will show the existence
of many
{\it doubly degenerate} groups, that is, surface groups whose
limit sets constitute the entire sphere.  We will construct
such groups as limits of quasi-Fuchsian groups where the conformal
structures on the two components of the domain of discontinuity
go to infinity in different directions.

\begin{theorem}[Double limit theorem]\label{double limit theorem}
Let $\mu$ and $\mu'$ be any two laminations in $\ML_0(S)$
satisfying the condition that for all $\nu \in \ML_0(S)$,
$i(\mu, \nu) + i(\mu', \nu) > 0$.

Then for any sequence $\{(g_i, h_i)\}$ in $\teich(S) \cross \teich(S)$
converging to $(\mu, \mu')$ in $\bar{\teich}(S) \cross \bar{\teich}(S)$,
the sequence of quasi-Fuchsian groups
$$ \{ \qf(g_i, h_i)\}$$
has a subsequence whose associated manifolds $\hy^3 / \qf(g_i, h_i)$
converge algebraically to a point $N_\infty \in \hs)$.
\end{theorem}

\begin{proof}[Proof of \ref{double limit theorem}]  First we will prove the theorem in the case that
$\bar S$ has non-trivial boundary, since this lacks some of the complications.

In this case, we can choose a lamination $\lambda$ of $S$ which is the
$1$-skeleton of an ideal triangulation. Represent the projective classes of
$\mu$ and $\mu'$ by measured laminations of the same names.

By Theorem \ref{laminations compactify Teichmuller space}, there are sequences
of measured laminations $\mu_i \arrow \mu$ and ${\mu_i}' \arrow \mu'$ such that
$$ \length_{g_i} ( \mu_i) \arrow 0$$
and
$$ \length_{h_i}( \mu'_i) \arrow 0.$$

Let $N_i$ be the quotient manifold of hyperbolic space by $\qf(g_i, h_i)$. From
either \ref{poincare length bounds hyperbolic length} or \ref{bounded
distortion to infinity}, it follows that $\length_{N_i}(\mu_i) \arrow 0$ and
$\length_{N_i}(\mu'_i) \arrow 0$. Let $P_{\lambda, i}$ denote the
$\lambda$-pleated surface in $N_i$.
We deduce from \ref{efficiency of pleated surfaces} that the lengths of $\mu_i$
and $\mu'_i$ on $P{_\lambda,i}$ remain bounded as $i \arrow \infty$. From
\ref{binding confinement} it follows that the hyperbolic structures of
$P_{\lambda,i}$ remain in a bounded portion of Teichm\"uller space. Therefore,
generators of the group $qf(g_i, h_i)$ remain bounded, since they move a base
point of $\tilde P{_\lambda,i}$ not further in $N_i$ than in $P_{\lambda,i}$. 
The theorem for the case that $\boundary \bar S$ is not empty follows.

\medskip

To prove the theorem when $S$ is a closed surface, we will show that for any
infinite sequence $\{N_i\}$ in $\hs$, we can always choose the lamination
$\lambda$ so that for an infinite subsequence the length of closed leaves of
$\lambda$ remains bounded above zero.

Since $\lambda$ can be chosen to have only one closed leaf, which can be any
simple geodesic on $S$, we only need to show that not all lengths of simple
geodesics tend to zero.

\begin{proposition}[Four not all small]\label{four not all small}
Suppose that $A$ and $B$ are isometries of $\hy^3$ which generate a discrete,
non-elementary group.  Then one of the four isometries $A$, $B$, $AB$ or
$AB\inverse$ is not conjugate to an isometry near the identity.
\end{proposition}

This can be proven by looking at $SL(2, \complexes)$, which is the double cover
of the group of isometries of $\hy^3$ and studying traces.  It is hard to keep
track of the geometry, however, in studying traces.  Instead, we will prove the
proposition using pleated surfaces.

\begin{proof}[Proof of \ref{four not all small}]
We can map the universal cover of the three-punctured sphere $M$ into $\hy^3$,
to represent the homorphism $\fund M \arrow G(A, B)$. The three-punctured
sphere has an ideal triangulation with two triangles, such that the three
corners of each triangle go to different punctures. Choose an orientation
for the axes of $A$, $B$ and $C = AB$ if these transformations are not
parabolic. We can define an equivariant map from the universal cover of this
triangulation to $\hy^3$ by mapping each corner of a triangle to the positive
endpoint of the corresponding axis in $\hy^3$, or to the parabolic fixed point
if it has no axis. The three corners of any triangle map to three distinct
points, for otherwise the group would be reducible, hence elementary. This
gives us a pleated surface, homeomorphic to a three-punctured sphere (but
geometrically, more likely the interior of a pair of pants), in the quotient
manifold or quotient orbifold.

The group can now be conveniently described in terms of the local geometry of
the pleated surface.  To each edge is associated a complex number, which may be
thought of as the edge invariant for the tetrahedron formed by the four
vertices of the two triangles which join at that edge.  An edge invariant of
$1$ means that the one triangle is obtained from the other by a $180\degrees$
rotation about the edge; an edge invariant of $-1$ means that the two triangles
are the same.  Denote the edge invariant for the edge going from boundary
component $A$ to boundary component $B$ by $c$, and so on.

Let $\alpha$, $\beta$ and $\gamma$ be the complex translation invariants for
$A$, $B$ and $C$: thus, $\log|\alpha|$ is the translation distance of $A$
(along its axis), and $\arg(\alpha$) is the rotation angle, etc.  Then, by
arranging an axis so its positive endpoint is at infinity in the upper
half-space model, it is easy to see that
$$ b c = \alpha, \  c a = \beta, \ \hbox{\rm and}\ a b = \gamma$$.

If $A$, $B$ and $C$ are conjugate to isometries near the identity, then
$\alpha$, $\beta$ and $\gamma$ are near 1.  There are algebraically two
possibilities for $a$, $b$ and $c$: they can all be near $1$, or all be near
$-1$.  But if they were all near $-1$, then the two triangles across any edge
would nearly coincide, and a coordinate system in $\hy^3$ could be chosen so
that $A$, $B$, and $C$ were simultaneously near the identity; this is
impossible for a discrete, non-elementary group.  Therefore, if $\alpha$,
$\beta$, and $\gamma$ are near $1$, the edge invariants are near $1$.  This
means that the triangles are nearly flat across the edges, so the group action
is near that of a standard three-punctured sphere.  But on a hyperbolic
three-punctured sphere, $A B\inverse$ is not near the identity.
\end{proof}

\begin{corollary}[Four curves not all short]\label{four curves not all short}
There is an $\epsilon$ such that for any $N \in \hs$, and any two simple closed
curves $\alpha$ and $\beta$ on $S$ such that $i(\alpha, \beta) = 1$, then the
length in $N$ of $\alpha$, $\beta$ or one of the two simple curves obtained by
cutting and pasting $\alpha$ and $\beta$ at their intersection point is greater
than $\epsilon$.
\end{corollary}
\begin{proof}[Proof of \ref{four curves not all short}]
In the fundamental group of $S$ based at the point of intersection of $\alpha$
and $\beta$, the homotopy classes of the four simple curves have the form $A$,
$B$, $AB$ and $AB\inverse$.  

We claim that if a simple curve on $S$ has a short length in $N$, then it also
has a small angle of rotation, so that its complex translation invariant is
near $1$.  To establish the claim, construct a pleated surface containing a
short closed geodesic $g$.  From the thick-thin decomposition of $N$, we know
that $g$ has a regular neighborhood in $N$ whose diameter goes to infinity as
the length of $g$ goes to zero. If $g$ has rotation angle $\theta$, then the
shortest curve homotopic to $g$ on the boundary of the solid torus of radius
$R$ has length at least $\theta \sinh R$.  In the covering of $N$ with
fundamental group generated by $g$, $\tilde S$ is divided by $g$ into two
semi-infinite cylinders, one of which intersects this torus in a cycle
homologous to $g$. If $g$ has a canonical regular neighborhood of radius $R$,
the area of the intersection of the pleated surface with this neighborhood is
therefore greater than $\int_0^R \theta \sinh t \ dt = \theta(\cosh R - 1)$. 
Since the area of the pleated surface is bounded, if $R$ is large, then
$\theta$ must be small.

Applying proposition \ref{four not all small} to the four elements of $\fund
S$, we conclude that they are not all conjugate to isometries near the
identity, and therefore their lengths in $N$ are not all small.
\end{proof}

The proof of theorem \ref{double limit theorem} can now be completed in
the case of a closed surface.  Suppose we are given a sequence $\{ ( g_i, h_i
)\}\in  \teich(S) \cross \teich(S)$ whose limit is $(\mu, \mu')$ as in the
statement of the theorem. We choose two curves $\alpha$ and $\beta$ on $S$
intersecting in one point, and an infinite subsequence such that one of the
four curves from the corollary, call it $\lambda_0$, always has length greater
than $\epsilon$ in $\qf(g_i, h_i)$. Extend $\lambda_0$ to a maximal lamination
$\lambda$ by adding a finite number of leaves spiraling at both ends around
$\lambda_0$, and consider the pleated surface $P_{\lambda, i}$ in the quotient
manifold $N_i$.

As before, there are measured laminations $\mu_i \arrow \mu$ and $\mu'_i \arrow
\mu'$ such that $\length_{N_i}(\mu_i) \arrow 0$ and $\length_{N_i}(\mu'_i)
\arrow 0$. From \ref{efficiency of pleated surfaces} we deduce that the
length of $\mu_i$ and $\mu'_i$ on $P_{\lambda,i}$ remain bounded as $i \arrow
\infty$. From \ref{binding confinement}, it follows that the hyperbolic
structures of $P_{\lambda,i}$ remain in a bounded portion of Teichm\"uller
space, and hence that some subsequence of the $N_i$ converges algebraically to
a limit $N_\infty \in \hs$.

\end{proof}

For easier reference, we collect in one place the main ingredients of this
proof:

\begin{theorem}[Surfaces approximate three-manifolds]\label{surfaces approximate three-manifolds}
For every surface $S$ equipped for reference with a complete hyperbolic
structure $h$ of finite area, there is a constant $C$ such that for any
hyperbolic three-manifold $N \in \hs$, there is a complete hyperbolic structure
of finite area $g$ on $S$ and an isometric embedding $f_\lambda$ of $S$ as a
pleated surface in $N$, such that for all laminations $\mu \in \ML_0(S)$,
$$ \llength_N(\mu) \le \length_g(\mu) \le \llength_N(\mu)
+ C \length_h(\mu).$$
\end{theorem}
\begin{proof}[Proof of \ref{surfaces approximate three-manifolds}]
Apply \ref{four curves not all short} to conclude that there are four specific
laminations $\lambda_i \ [i=1,\dots,4]$ such that in any element $N \in \hs$,
at least one of the four satisfies the hypotheses of \ref{efficiency of pleated
surfaces}.  The alternation number $a(\lambda_i , \mu)$ is bounded by some
multiple of $\length_h(\mu)$, for each $i$; hence, the inequality.
\end{proof}

Of course, an equivalent statement would be obtained if $\length_h(\mu)$ were
replaced by wordlength with respect to some fixed set of generators.


\section{Hyperbolic structures for mapping tori}
\label{mapping tori}

We will now prove the main theorem,
\ref{mapping torus hyperbolic}, which we restate:

\begin{maintheorem}
Let $M^3$ be a compact $3$-manifold (possibly with boundary) which
fibers over
$S^1$, and whose fiber is a compact surface of negative Euler characteristic.

Then the interior of $M$ either

\begin{description}
\item[(i)] has a complete $\hy^2 \cross \reals$ structure of finite volume,
and can be described as a Seifert fibration over some hyperbolic $2$-orbifold,
\item[(ii)] contains an embedded incompressible
torus not isotopic to a boundary component,
and splits along this torus into two simpler three-manifolds, or
\item[(iii)] (generic case) has a complete hyperbolic structure of finite volume.
\end{description}

Cases (i) and (ii) are not mutually exclusive, but (iii) excludes
the other two cases.
\end{maintheorem}

\begin{proof}[Proof of \ref{mapping torus hyperbolic}]
Let $M^3$ be a 3-manifold which fibers over the circle, and let $\phi$ be its
monodromy. By \ref{non-pseudo-Anosov mapping torus}, either $M_\phi$ has
a complete, finite volume $\hy^2 \cross \reals$ structure, is torus-reducible,
or $\phi$ is pseudo-Anosov. What remains to be proven is that if $\phi$ is
pseudo-Anosov, then the interior of $M_\phi$ admits a complete hyperbolic
structure of finite volume.

Choose an arbitrary element $g \in \teich(S)$, and consider the
sequence of quasi-Fuchsian groups $\Gamma_i = qf(\phi ^{-n} g, \phi^n g)$.
By \ref{classification of surface homeomorphisms}, the sequence
$\{ (\phi ^{-n} g, \phi^n g) \}$
tends to the limit $(\mu_s, \mu_u)$
in $\overline {\teich(S)} \cross \overline {\teich(S)}$, where $\mu_s$ and
$\mu_u$ are the stable and unstable laminations of $\phi$.
This pair of laminations satisfies the hypothesis of
\ref{double limit theorem},  so there is a subsequence of the $\Gamma_i$
which converges to a group $\Gamma_\infty$.

Note that all the hyperbolic surfaces $\phi^n(g)$ are isometric: they differ
only by their parametrization, or ``marking''. Furthermore, $\phi$ acts
geometrically as the same map on each of the surfaces, since it commutes with
any iterate of itself. In particular, $\phi$ has a fixed quasi-isometric
constant on each of the hyperbolic surfaces $\phi^n(g)$.

If we lift these actions of $\phi$ to the two components of
the domain of discontinuity for $\Gamma_i$, we obtain a homeomorphism $h_i$
of the sphere, equivariant twisted by $\phi$
with respect to $\Gamma_i$.  (That is, $h ( \gamma(x)) = \phi(\gamma) (h(x))$).
This homeomorphism is quasiconformal (by general principles), and 
its quasiconformal
constant is bounded uniformly over $i$, since it is uniformly bounded
on the domain of discontinuity, which has full measure: in fact, the
best constant is a function only of the pseudo-Anosov constant.

The groups $\Gamma_i$ are really only defined up to conjugacy by
Moebius transformations.  To pin them down as actual groups,
pick three non-commuting hyperbolic elements, and arrange that
the attracting fixed points of these three elements are at $0$,
$1$ and $\infty$ on $\sinfty$.  The actual groups now converge,
and the points $\phi(0)$, $\phi(1)$, and $\phi(\infty)$ also converge
to three distinct fixed points
of the limits of the images of the three elements by $\phi$.
(In the limit, the elements conceivably
might not be hyperbolic, but that is
not a difficulty).

The sequence of maps $h_i$ is uniformly quasiconformal, and converges
on the three points $0$, $1$ and $\infty$ to three distinct limits.
Therefore, the
sequence is equicontinuous, and there is a convergent subsequence.
The limit is a quasiconformal homeomorphism $h_\infty$
which is equivariant
twisted by $\phi$
with respect to $\Gamma_\infty$.

We claim that the limit set of $\Gamma_\infty$ is all of $\sinfty$.
One proof of that fact goes as follows.
By the Ahlfors finite area theorem, the quotient of $D_{\Gamma_\infty}$
by $\Gamma_\infty$ is conformally equivalent to
a (not necessarily connected) hyperbolic surface $A$ of finite area.
The homeomorphism $h_\infty$ descends to a homeomorphism of this
quotient surface.  If there were a component of $A$ which had a cusp,
there would be a finite number of cusps, which were permuted by $h_\infty$;
this would contradict the fact that $\phi$ is pseudo-Anosov.
The other possibility for a component of $A$ is that it is a finite
sheeted covering space of $S$.  If such existed, then the conformal
structure induced on a finite sheeted cover of one component of the
domain of discontinuity would converge; that is absurd.

Now we can apply theorem \ref{quasiconformal deformations},
to conclude that the quasiconformal homeomorphism $h_\infty$ actually
induces a conformal automorphism of $\Gamma_\infty$ to itself.
This defines a representation of $\fund N$ in the group of conformal
homeomorphisms of $\sinfty$, where $\fund S$ acts as $\Gamma_\infty$
and $\phi$ acts as $h_\infty$.  We claim that the group is discrete
and faithful.  To see this, consider the construction in terms of
quotient manifolds.  First form the quotient $N_\infty$ of
$\hy^3$ by $\Gamma_\infty$; it is homotopy equivalent to the
surface $S$.  The conformal map $h_\infty$ extends to an
isometry of hyperbolic space, which descends to an action as an
isometry on
$N_\infty$.  But the group of isometries of any hyperbolic manifold
with non-elementary fundamental group is easily seen to be discrete,
so therefore the action of $\fund N$ is discrete.

By a theorem of Stallings \cite{St}, the quotient hyperbolic manifold
is in fact homeomorphic to $N$.
\end{proof}


\section{On limits and limiting behavior of surface groups}
\label{surface groups}

We have seen that if a sequence $\{(g_i, h_i)\}$
in $\teich(S) \cross \teich(S)$
converges to a point at infinity $(\mu, \mu')$ far from the diagonal,
then the corresponding sequence of quasi-Fuchisan groups has at
least a subsequence that converges algebraically.
On the other hand, the quasi-Fuchsian groups corresponding to
a diagonal sequence $(g_i, g_i)$ which tends to infinity can
never converge.

Let us reformulate this by describing a compactification
of $\hs$ which adjoins abstract limits to all non-convergent
sequences.

A point $N \in \hs$ gives rise to a continous length
function
$$ \length_N = l_N: \ML_0 \arrow \reals$$
which we can think of as defining a point
$l^*(N)$ in the vector space $\reals^{\ML_0}$.
To consider how ratios of lengths behave in hyperbolic manifolds,
we project
$\reals^{\ML_0}$ into its associated projective space
\def\prm{\proj (\reals^{\ML_0})} $\prm$, and write
$$\proj l^* : \hs \arrow \prm$$
for the composition.
Let $\prm$ have the quotient topology of the compact-open topology
on $\reals^{\ML_0}$.

\begin{theorem}[Closure compact]\label{closure compact}
The closure of the image of $\proj l^*$ is compact.  The added
\def\OImage{\mathop{\overline{\hbox{\rm Image}}}\nolimits}
\def\Image{\mathop{\hbox{\rm Image}}\nolimits}
set $\OImage  \proj l^* -
\Image \proj l^*$
is the same as for Teichm\"uller space:
it is $\PL_0$, as it is embedded in $\prm$ via geometric intersection
number.
\end{theorem}

In other words, for any sequence of hyperbolic three-manifolds $N_i$
in $\hs$, either there is an algebraically convergent subsequence,
or there is a subsequence such that the ratios of lengths of simple
closed curves (and laminations) all converge, and there is
a sequence of Fuchsian groups with the same limiting ratios!

\begin{proof}[Proof of \ref{closure compact}]
The proof is a corollary of what we have already done.  Apply
Theorem \ref{surfaces approximate three-manifolds} to find
a uniformly efficient pleated surface in each member
of the sequence of
quotient manifolds $N_i$, that is, satisfying
$$ \length_{N_i} (\mu) \le \length_{P_\lambda}(\mu) \le 
\length_{N_i} (\mu) + C a(\lambda, \mu) .$$

If some of the lengths tend to infinity, the deviation bounded by
$C a(\lambda, \mu)$ becomes irrelevant in projective sapce.
There is either a subsequence that converges algebraically,
or, a subsequence such that
the  sequence of hyperbolic structures on $S$ defined by $P_\lambda$
converge to a point in $\PL_0$.  The sequence of
length functions defined by $N_i$ then converge to the image point
in $\prm$.
\end{proof}
\bigskip

The limiting behavior of the length functions in $\prm$ is only
part of the story.
Sometimes one needs to understand infinite nonconvergent sequences
of Fuchsian or quasi-Fuchsian groups in a more refined way.

\begin{theorem}[Converge on subsurface]\label{converge on subsurface}
Let $\{N_i\}$ be any sequence of elements of $\hs$
Then there exists a subsequence $N_{\{c(j)\}}$ and a possibly
empty, possibly disconnected subsurface $S'$ of $S$ with incompressible
boundary, such that 

(a)  For each component $S'_i$ of $S'$, the sequence of representations
of $\fund {S'_i}$ in  $\Isom \hy^3$ coming from $N_{c(j)}$ converges
up to conjugacy.

(b)  If $\Gamma$ is any nontrivial subgroup of $\fund S$ such that
for some subsequence of indices $\{c(d(j))\}$ of $\{c(j)\}$,
its sequence of representations
converges up to conjugacy,
then $\Gamma$ is conjugate to a subgroup of $\fund {S'_i}$ for
some $i$.
\end{theorem}

The subsequence $\{c(j)\}$ is a {\it sequence of
maximal convergence}, and
the subsurface $S'$ is the
{\it surface of maximal convergence} associated to that subsequence.
There may be various subsequences of maximal convergence contained
in a given sequence, associated with different subsurfaces.

\begin{proof}[Proof of \ref{converge on subsurface}]
This theorem is reduced to the Fuchsian case, by an application
of 
\ref{efficiency of pleated surfaces}
and \ref{four curves not all short}.

To prove the theorem in the Fuchsian case, that is, for a sequence
of hyperbolic structures on $S$, we proceed as follows.

To any collection $A$ of homotopy classes of closed curves on a surface is
associated a subsurface $\ls(A)$, unique up to isotopy, which is
the least surface with incompressible boundary that contains
$A$.  It is characterized by the properties
that
\begin{description}
\item[(a)] Every $\alpha \in A$ is isotopic to a curve in $\ls(A)$
\item[(b)] No component of $\ls(A)$ is an annulus isotopic into another
component,
\item[(c)] For any other subsurface $R$ of $S$ satisfying (a), $\ls(A)$ is
isotopic to a subsurface of $R$, and
\item[(d)] $\boundary \ls(A)$ is incompressible.
\end{description}

Here is one construction for $\ls(A)$.  Choose a hyperbolic structure on $S$,
and represent each element $\alpha \in A$ by its geodesic, or a small horocycle
in the case that the element is parabolic. Now fill in each component of the
complement which is a disk. Also, fill in, one by one, any region of the
complement which is an annulus one of whose boundary components is a
component of the set so far constructed. Thicken the result, so that it is a
subsurface.  It is easy to verify that it has the desired properties.

We claim that if we have a sequence of hyperbolic structures of on a surface
$S$ and if the lengths of all homotopy classes of curves in a collection $A$
remain bounded, then the sequence of representations of the fundamental group
of any component of the surface $\ls(A)$ converges.

The basic observation is that two simple closed geodesics $\alpha$ and $\beta$
intersect in a combinatorial pattern which depends only on their homotopy
classes.  If the lengths of $\alpha$ and $\beta$ remain bounded in a sequence
of hyperbolic structures, then $\alpha \cup \beta$ forms a $1$-skeleton for the
surface $\ls(\{\alpha, \beta\})$, so the representations of the fundamental
group of that entire subsurface remain bounded, and a subsequence can be
chosen so that the sequence of representations converge.

Unfortunately, the combinatorial pattern of the intersection of three or more
simple closed geodesics is not invariant: it depends on the hyperbolic metric.
However, an inductive argument can be used to show that if $A$ is any
collection of simple closed curves on a surface whose lengths all remain
bounded, then the sequence of representations of the fundamental group of any
component of $\ls(A)$ converges.  In fact,  if $S_1$ is the least subsurface
for some finite subcollection $A_1$, and if $\alpha$ is an additional curve,
then it is possible to find a collection $A_0$ of disjoint simple curves in
$\ls(A_1)$ such that $\ls(\{A_0\} \cup \alpha) = \ls(\{A_1\} \cup \alpha)$. 
Note that $A_0$ need not be a subset of $A$. The geodesics in the $\{A_0\} \cup
\alpha$ provide a bounded 1-skeleton for this surface.

This takes care of collections of simple curves, but not of general curves. 
For any finite collection $F$ of homotopy classes of loops on a surface, there
is a finite-sheeted covering of the surface where the induced covers of each
loop in the collection are embedded. This is a consequence of a theorem of Peter
Scott \cite{Scott:LERF},
although this particular case is easier than the general
theorem. This implies that the sequence of representations of a subgroup of
finite index in the fundamental group of any component of $\ls(F)$ remains
bounded if the lengths of elements of $F$ remain bounded, so a subsequence can
be chosen so that these representations converge.  But hyperbolic
transformations of $\hy^2$ are uniquely divisible, so the representation of a
subgroup of finite index determines the representation of the full group.
\end{proof}

Theorem \ref{converge on subsurface} is related to
proposition \ref{binding confinement}, and they can be combined into
a single statement which generalizes
\ref{double limit theorem}.  If $\mu$ and $\nu$ are geodesic
laminations, let us say that $(\mu, \nu)$ binds a subsurface $R \subset S$
if
$\boundary R$ is incompressible in $S$, each component of $R$
has negative Euler characteristic, and $(\mu, \nu)$ binds $R$
equipped with any complete hyperbolic metric of finite area.

In this situation,
if $\mu \subset \mu_1$ and $\nu \subset \nu_1$, we will
also say that $(\mu_1, \nu_1)$ also bind $S$.  The purpose of
this is to allow for components of one lamination which do
not intersect the other.

\begin{theorem}[Laminations bind subsurface]\label{laminations bind subsurface}
If $N_i \in \hs$ is a sequence of hyperbolic manifolds
and $\mu _i\arrow \mu, \nu_i \arrow  \nu \in \ML_0$ are
sequences of measured laminations such that
$ \length_{N_i} (\mu_i) and \length_{N_i} (\nu_i)$ remain bounded,
then if $(\mu, \nu)$ bind a subsurface $R \subset S$, there is
a subsequence such that the sequence of representations of
the fundamental group of $R$ converges.
\end{theorem}

\begin{proof}[Proof of \ref{laminations bind subsurface}]
As before, we need only prove the theorem in the Fuchsian case, since
under the hypotheses
we can construct a pleated surface on which both $\mu_i$ and
$\nu_i$ remain bounded in length.

We will exploit the fact that the combinatorial pattern of
intersections of $\mu$ with $\nu$ does not depend on a hyperbolic
structure on $S$.
We will phrase the proof for the case that $\mu$ and $\nu$ have
no measure concentrated on individual leaves.  The general case
works in the same way, but it is awkward to word the proof to apply
to all cases at once.

Let $\mu_0$ and $\nu_0$ be components of
$\mu$ and $\nu$, respectively, which actually intersect each other.
It suffices to construct a subsequence which
works for the surface which $R_0$ $(\mu_0, \nu_0)$
binds, in view of \ref{converge on subsurface}.

Let $x$ be a point in $\mu_0 \cap \nu_0$.  A small open and closed
neighborhood $V$ of
$x$ in the intersection can be formed which has the structure
of the product of the local leaf space for $\mu$ with the
local leaf space for $\nu$.  $V$ is contained
in a quadrilateral $Q$ on $S$ whose sides are leaves of $\mu$ and $\nu$.

Every leaf in $\mu_0$ and in $\nu_0$ eventually intersects $Q$, in each
direction.
The leaves of $\mu_0$ and of $\nu_0$ minus $Q$ therefore
group into a finite
number of bands which remain roughly parallel.  Each of the bands
has a positive transverse measure.  We will work with extended
bands, where both ends are extended across $Q$.

Consider the set of hyperbolic structures on $S$ for which the
length of $\mu_0$ and the length of $\nu_0$ both are less than a
constant $C$.  Then the average length of a leaf in one of the extended
bands is less than $2 C$ divided by the transverse measure of
the band  --- a uniform constant.   Select one leaf from
each extended band whose length does not exceed the average length.

The union of the selected leaves forms a 1-skeleton for
the subsurface $R_0$.  The combinatorial type outside of $Q$
is fixed; $Q$ is simply connected, and
inside $Q$, any two points can be connected
by a path of bounded length. 
This proves that the sequence
of representations of $\fund {R_0}$ has a convergent subsequence,
under the hypothesis
that the lengths of $\mu$ and $\nu$ remain bounded.

This proof generalizes without difficulty
to the case where $\mu_i \arrow \mu$ and
$\nu_i \arrow \nu$, with the lengths of $\mu_i$ and $\nu_i$ bounded
by a constant $C$.  To do this, consider a neighborhood in
$S$ slightly bigger than $Q$, but still intersecting the same leaves
of $\mu$.  For each band of leaves of $\mu$ there is a nearby band
of nearly parallel leaves of $\mu_i$, whose total transverse measure cannot
decrease suddenly. There may also be new bands of small measure hitting
$Q$,
but we don't need to use them.  As before, we can construct a
1-skeleton for $R_0$, made up from leaves
of $\mu_i$ and $\nu_i$.

The proof is completed by letting $\mu_0$ and $\nu_0$ range over
all components of $\mu$ and $\nu$ and
applying \ref{converge on subsurface}. (An alternative would
be to use several quadrilaterals).
\end{proof}
\bigskip

We give some more examples illustrating
how the double limit theorem can be applied to
construct Kleinian groups.

\begin{theorem}[Accidental parabolics]\label{accidental parabolics}
Let $\alpha$ and $\beta$ be collections of disjoint
simple closed geodesics on $S$.  If
$\alpha \intersect \beta = \emptyset$,
then there is a sequence of elements of $\QF(S)$ converging
to an element of 
$$\AH(S \cross I, \boundary S \cross I \cup
\Nbhd(\alpha) \cross \{0\} \cup \Nbhd(\beta) \cross \{1\}) \subset \hs.$$
\end{theorem}

\begin{proof}[Proof of \ref{accidental parabolics}]
Construct a sequence $\{g_i\}$ of hyperbolic structures
on $S$ such that it is a sequence of maximal convergence associated
with the surface $A = S - \Nbhd(\alpha)$.

Construct a similar sequence $\{h_i\}$ for $\beta$, converging
on the surface $B = S - \Nbhd(\alpha)$.

Now apply \ref{converge on subsurface} to
$\qf(g_i, h_i)$. 
The subsurface of maximal convergence must contain both $A$
and $B$, so it must be all of $S$.  Thus there is a subsequence
converging to a limiting group $\Gamma$.

Each of the curves in $\alpha$ and $\beta$ is parabolic, so the
manifold is in $\AH(S \cross I, \boundary S \cross I \cup
\Nbhd(\alpha) \cross \{0\} \cup \Nbhd(\beta) \cross \{1\})$.
\end{proof}

Here is a variation of this construction, where
$\alpha$ and $\beta$ are laminations which do not necessarily
bind $S$, or even a subsurface of $S$.
The statement is
more complicated because not every measured lamination whose
length in $N_i$ remains bounded is contained in the subsurface
of convergence.

\begin{proposition}[Convergence on non-binding laminations]\label{convergence on non-binding laminations}
Let $\lambda_i \ [i=1,2]$ be measured laminations on $S$.
For any complementary region of $\lambda_i$ and any ideal polygonal curve
on its boundary which is homotopically non-trivial, suppose that
there is a closed leaf of $\lambda_i$ in its homotopy class.

Then if $\lambda_1$ has no leaves in common with $\lambda_2$,
there is a sequence of elements in $\QF(S)$ converging
to an element of $N \in \hs$ in which $\lambda_1$ and
$\lambda_2$ are unrealizable, and
$$ \llength_N (\lambda_1) = \llength_N(\lambda_2) = 0$$.
\end{proposition}

\begin{proof}[Proof of \ref{convergence on non-binding laminations}]
Let $A_i \subset S$ be the possibly empty subsurface formed as
the union
of all components of $S - \lambda_i$ whose boundary consists of
closed geodesics.  Let ${g^i}_j \  [i=1,2]$ be a sequence
of maximal convergence
of hyperbolic structures on $S$ such that $A_i$ is its subsurface
of maximal convergence, and the length of $\lambda_i$ tends to
zero.

Form a subsequence of maximal convergence for the sequence of
quasi-Fuchsian groups $\qf({g^1}_j, {g^2}_j)$.

Let $\mu_i \subset \lambda_i$ be the union of components which
intersect $\lambda_j$ where $j \ne i$.
Then $(\mu_1, \mu_2)$ bind a certain subsurface $S_0$ of $S$.
The subsurface of maximal convergence must contain $A_1$, $A_2$, and $S_0$.
The only possibility is all of $S$.
\end{proof}


\section{Infinitely generated geometric limits}
\label{geometric limits}

In this section, we will describe a
construction for a different kind of limit
for surface groups, of a type first discovered
by Troels J\o rgensen ---
geometric limits
of sequences of quasi-Fuchsian manifolds ---
whose fundamental groups are not finitely generated.

One reason for studying geometric limits, as well as algebraic limits,
is that only the geometric limit captures the quasi-isometric information
of the approximating manifolds.  So far, there are three known invariants
which can distinguish Kleinian groups (at least those groups which
are not free products):  the topology, the conformal structure
on the quotient of the domain of discontinuity, and the ending
lamination (defined in \refin{Ch. 8 and 9}{GT3M}, and now known
to exist in a general setting, by the main result of \cite{Bon}).
If we can precisely analyze the quasi-isometric structure
of surface groups, we can probably resolve the question of whether
these three invariants are sufficient to determine a group.
The most difficulty seems to come from surface groups
which have arbitrarily short geodesics --- that is precisely
the situation which we will deal with in this section.

Let $S$ be a hyperbolic surface of finite area.
We will study finite and infinite
collections $C$ of disjoint
non-trivial
simple curves on levels in $S \cross \reals$,
$C = \set{\gamma_i = \alpha_i \cross \set t_i } $, 
up to isotopy through similar systems of curves.  We assume that
no two curves of $C$ are isotopic.  We also assume that $C$ is closed,
so that the set of levels involved is discrete.  Given such a collection
$C$, define $M_C = S \cross \reals - \bigcup C$.

There is a partial ordering on $C$, defined by
$ \gamma_i \prec \gamma_j$ if $t_i < t_j$, and if further, this inequality
remains true after an arbitrary isotopy.  Clearly, if $i(\alpha_i, \alpha_j)$
is not zero, then the two curves $\gamma_i$ and $\gamma_j$ are comparable,
and ordered according to $t_i$ and $t_j$.  Furthermore, if
there is a finite sequence $\set c(k)$ such that
$i( \alpha(c(k)), \alpha(c(k+1)) > 0$, and $t_{c(k)} < t_{c(k+1)}$,
then the elements of the sequence $\gamma_{c(k)}$ are in ascending
linear order.

Conversely, if there is no such chain connecting
two elements of $C$ (in either order), then they are incomparable.
Any set of mutually incomparable elements is isotopic to
a single level.  The maximum size of such a set is bounded
as a function of the topology of $S$.
On the other hand, it might happen that some element
of $C$ is incomparable to an infinite collection of other elements.

If $\gamma_1 \prec \gamma_2$, and if there is no intervening element
$\gamma_3$ such that $\gamma_1 \prec \gamma_3 \prec \gamma_2$, then
$\gamma_2$ is a {\it successor} of $\gamma_1$.  The partial ordering
is generated by the successor relation, as its transitive closure.
The successors of an element are incomparable, so the number is
bounded as a function of the topology of the surface.
The partial ordering can be represented by a directed acyclic graph,
with an edge joining each element to each of its successors.

A {\it cut} of $C$ is a partition of $C$ into two
parts $L$ and $U$, such that
\begin{description}
\item[(1)]no element of $L$ is greater than any element of $U$,
\item[(2)]every linearly ordered subset of $L$ contains a maximum, and
\item[(3)]every linearly ordered subset of $U$ contains a minimum.
\end{description}

The cuts themselves have a partial order, with
$(L,U) \prec (L', U')$ if L is a proper subset of $L'$.

\begin{proposition}[Surfaces cut]\label{surfaces cut}
For every cut $K = (L, U)$ of $C$, there is an isotopy of $C$
so that $L$ consists of the curves below $S \cross \set 0 $,
and $U$ consists of the curves above $S \cross \set 0$ .
There is associated to
$K$ an isotopy class $S_K$ of surfaces
in $M_C$ which separates $L$ from $U$,
and is isotopic in $S \cross \reals$ to $S$.  Conversely, every isotopy
class of embeddings of $S$ which project
homeomorphically under the map $M_C \arrow S$
defines a cut.
\end{proposition}

\begin{proof}[Proof of \ref{surfaces cut}]
There can originally be only be a finite number of elements of $L$ above $S
\cross 0$: if there were an infinite number, then there would be an infinite
sequence which is physically ascending, therefore non-descending in the partial
order.  But the property that the maximal size of pairwise incomparable set is
bounded implies that every infinite non-descending sequence contains an
infinite ascending subsequence.

Similarly, there are only a finite number of elements of $U$ below $S \cross 0$.

We can repeatedly move the physically lowest element of $L$ which is above $S
\cross 0$, until there are none left, because the elements of $U$ cannot
interfere.  Then we can do the same process for $U$, proving the proposition.

After the isotopy of $C$, $S \cross 0$ serves as the surface $S_K$.  Because
all isotopies of the system of curves can be arranged so that they involve only
vertical motion, the isotopy class of $S_K$ is well-defined.
\end{proof}

Note that the partial ordering on cuts corresponds to
the partial ordering on surfaces: if $A$ and $B$
are two surfaces, $A \prec B$ if $B$ can be isotoped below
$A$, but $A$ and $B$ are not isotopic.

\begin{theorem}[Drill holes]\label{drill holes}
Let $C$ be a collection of curves in $S \cross \reals$ as above.
Suppose that $C$ has the property that
for every non-trivial simple closed
curve $\beta$ on $S$, if $C$ is unbounded above then
the set of $\alpha \cross t \in C$ such that $i(\alpha, \beta) > 0$
is unbounded above, and if $C$ is unbounded below then
the set of $\alpha \cross t \in C$ such that $i(\alpha, \beta) > 0$
is unbounded below.

Then there is a sequence of hyperbolic
manifolds in $\QF(S)$ whose geometric limit is a manifold $N$
homeomorphic to $M_C = S \cross \reals - \bigcup C$.
In particular, if $C$ is infinite, $\fund N$ is not finitely generated.
\end{theorem}

\begin{remarks}
A similar analysis would seem to work if the hypothesis
on $C$ were dropped, but the picture is somewhat different.  In the
present situation, the new cusps in the geometric limit are all
$\integers \cross \integers$-cusps; in the more general situation, there
would also be new $\integers$-cusps.

The proof of this theorem is probably of more interest than the statement; it
gives techniques for finding when there are short geodesics in a hyperbolic
manifold with the fundamental group of a surface. Eventually, it would be nice
to have a method so that, given two conformal structures at infinity, and given
any measured lamination on the surface, we can say, to within some bounded
factor, what is the length of the lamination in the resulting quasi-Fuchsian
manifold.  A prerequisite is probably a better understanding of the geometry of
the mapping class group of a surface, and of the various geometries of
Teichm\"uller space.

There is, in fact, a cheap proof of the current theorem at least for many
families of curves $C$. Here, in outline, is how it goes. First, construct by
some device a hyperbolic manifold homeomorphic to $M_C$.  For example, suppose
that $C$ can be isotoped to the form that on each level where a curve in $C$
occurs, there are enough other elements to fill up the surface. Then we can
construct a geometrically finite surface group corresponding to the region
between any two such levels,  such that the curves above and below are
parabolic, by citing \ref{accidental parabolics}.  All the stabilizers of
domains of discontinuity are three-punctured spheres, which implies that the
convex core of these surface groups has totally geodesic boundary, each
component of which is isometric to a standard three-punctured sphere. Glue
together three-punctured spheres, as appropriate, to obtain the limit.

Once a limit manifold is constructed, it is possible to work backwards, at
least in a case such as the one described, with the aid of the theory of
hyperbolic Dehn filling.  Throw away all the pieces which have been assembled,
except for a finite number in the middle; then hyperbolic Dehn filling can be
performed, to obtain a quasi-Fuchsian group.

Alternatively, one can use the general existence of hyperbolic structures on
Haken manifolds to obtain the manifolds with finitely many stages (to be proven
in part IV of this series), and use the deformation theory for three-manifolds
which have cylinders, (to be presented in part III of this series) to conclude
that as more and more stages are added, there is a geometric limit homeomorphic
to $M_C$.  Each of the finite stages is a limit of quasi-Fuchsian groups, using
hyperbolic Dehn filling, and one can diagonalize to show that $M_C$ itself is a
geometric limit.
\end{remarks}

Rather than taking an external approach such as this, we will give an internal
proof: we will begin with the quasi-Fuchsian groups, and find enough ways to
analyze their quasi-isometric geometry to understand the nature of the limit.

\begin{proof}[Proof of \ref{drill holes}]  The construction, in its geometric form,
goes as follows.

Index the elements of $C$ 
by an interval of positive and negative integers including 0,
so that the the physical height of the curves is monotone in
the indices.

For any curve $\gamma$ in $S \cross \reals$ which lies on a level,
let $A_\gamma$ be an annular neighborhood of $\gamma$ which
lies on the same level.  If we slit open $S \cross \reals$
along $A_\gamma$, and then reglue the two sides of the
cut by a power of a Dehn twist, we obtain a new manifold
which is still homeomorphic to $S \cross \reals$.
However, there are two natural isotopy classes
of homeomorphisms to the original
space:  there is one homeomorphism which is
the identity below the level of $\gamma$, and another
which is the identity above that level.
This construction is a special case of Dehn surgery
along $\gamma$.

Suppose we perform such a construction, using various powers
$q(i)$
of Dehn twists (perhaps sometimes the $0$th power),
on all the elements of our set $C$.
We obtain a new manifold $M_q$.
For each cut $K$, there
is an isotopy class of homeomorphisms $h(q, K)$
with the
original model so that $S_K$ maps as the identity.
A pair of cuts $K$ and $K'$ determine an isotopy class of homeomorphisms
$\phi(q, K, K')$
from $S$ to itself, which give the comparison:
$ h(q, K') \compose \phi(q, K, K') \homotopic h(q, K)$.

We will arrange a sequence of quasi-Fuchsian manifolds
so that the conformal structures on the top and bottom
are constant {\it up to diffeomorphism}.
Choose a conformal structure, and think of it as being
on $S \cross +\infty$ and $S \cross -\infty$.
We will choose the powers $q$ so that all but a finite
number are $0$.  Peform Dehn surgery to obtain $M_q$.
The homeomorphisms $h(q, K)$ all agree for $K$ sufficiently
high; define this as $h(q, +\infty)$.  Define $h(q, -\infty)$,
similarly, and also extend the definition
of $\phi(q, K, K')$ to allow $K$ or $K'$ to be $\pm \infty$.
If $K$ is an arbitrary cut,
then a quasi-Fuchsian group $\Gamma(q, K)$ is defined
by pushing the conformal structures at the two ends of
$M_q$ to $S \cross I$ via the homeomorphism $h(q, K)$:
that is,
$$\Gamma(q, K) = qf( \phi(q, +\infty, K)(g), \phi(q, -\infty, K)(g) )$$.

The groups for the various cuts are all isomorphic, of course.
Another way to think of this construction is that we have defined
a sequence of representations
$$\rho_q: \fund {M_C} \arrow \Isom ( \hy^3),$$
non-faithful with a quasi-Fuchsian image,
and the parametrization by $\fund S$ associated with the cut $K$
is induced by the embedding of $S$ in $M_C$ as $S_K$.

We will use sequences of powers $q_j(i)$ which are 0 outside a
finite interval of integers $[-n, m]$, and which are
large within that interval.  Given such a choice,
what do the conformal structures
at $+\infty$ look like from the perspective of a surface $S_K$?

\begin{proposition}[Big twists]\label{big twists}
Let $\alpha$ be a measured lamination which is an integer linear
combination of simple closed curves, and let $D$ be a diffeomorphism
which is the composition of either ${\tau_\gamma}^n$ or ${\tau_\gamma}^{-n}$
where $\gamma$ ranges over the components of $\alpha$ and $n$
is the weight of $\gamma$.

Write $\alpha = \sum  \alpha(i)$, where $\alpha(i)$  ranges over the
components of $\alpha$.
If $\lambda$ is any measured lamination,
then
$$ \lim_{q \arrow \infty} \frac {1}{q} D^q (\lambda) = \sum i(\lambda, \alpha(i)) \alpha(i)$$
in $\ML_0(S)$.

If $g$ is any hyperbolic structure on $S$, then
$$ \lim_{q \arrow \infty} \frac {1}{q} D^q (g) = \sum \length_g (\alpha_i) $$
in $\overline{\teich(S)}$
\end{proposition}

\begin{proof}[Proof of \ref{big twists}]
If $\lambda$ is a simple closed curve, then the result of $D^q$ applied
to $\lambda$ is to modify $\lambda$ so that it winds $n(i) q $ times
around the support of $\alpha(i)$, each time it intersects.  Clearly
the limit is as stated. 

The case of general $\lambda$ could be proven in
a similar way with the aid of train tracks, but it can also be 
logically derived,
as follows, from the case of simple closed curves.  Any measured lamination
is determined by its intersection numbers with simple closed curves.
To show that we have the right formula for a general lamination $\lambda$,
it suffices to show that the intersection numbers with simple closed
curves have the correct limits. 

For any measured lamination $\lambda$, write 
$A(\lambda) = \sum i(\lambda, \alpha(i)) \alpha(i)$, the expected
answer.  We have
$$i(\lambda, A(\lambda')) = i(A(\lambda), \lambda'),$$
and also 
$$ i(D^q (\lambda), \lambda') = 
i(\lambda, D^{-q} \lambda') .$$
Letting $\lambda'$ range over simple closed curves, these equations
imply that
$A(\lambda)$ has the correct limiting intersection with all simple
closed curves, so it is the correct limit.

Similarly
if $g$ is a hyperbolic structure,
define $A(g) = \sum \length_g(\alpha_i) \alpha_i$.  Again,
we have the symmetry $i(A(g), \lambda)$ = $\length_g(A(\lambda))$.
As before, this proves that $A(g)$ is the right answer for
the limit of $D^q(g)$.
\end{proof}

\begin{proposition}[Compose big twists]\label{compose big twists} 

Let $A = \set {\alpha(i)|i=1\dots n} $ be a sequence of simple closed curves
on $S$.  

If $\mu \in ML_0$ is a measured lamination such that each of
its components has positive intersection number with
at least one element of $A$, then for sufficiently large positive or
negative integers
$\set {q(i)|i=1,\dots,n} $, the image of $\mu$ by the
composition
$$({\tau_{\alpha(1)}}^{q(1)} \compose {\tau_{\alpha(2)}}^{q(2)}
\compose \dots \compose {\tau_{\alpha(n)}}^{q(n)}) (\mu)$$
is near the simplex in $\PL_0(S)$ generated by those curves
$\alpha_i$ which do not intersect any earlier element of the
sequence.

Similarly, if $g$ is any element of $\teich(S)$, then
for sufficiently large integers
$\set {q(i)|i=1,\dots,n} $, the image of $g$ by the
composition
$$({\tau_{\alpha(1)}}^{q(1)} \compose {\tau_{\alpha(2)}}^{q(2)}
\compose \dots \compose {\tau_{\alpha(n)}}^{q(n)}) (\mu)$$
is near the simplex in
$\PL_0(S)$ (considered as the boundary of Teichm\"uller space)
generated by those curves
$\alpha_i$ which do not intersect any earlier element of the
sequence.
\end{proposition}

Note that the first elements of any such subsequences of maximal
length must have zero intersection number with all earlier elements in
$A$.
Therefore, the first elements of any pair of maximal subsequences
have zero intersection number, so that a linear combination
of them makes sense as a measured lamination, and gives
a well-defined projective class.  The set of all linear
combinations of them map to a simplex in $\PL_0(S)$.

\begin{proof}[Proof of \ref{compose big twists}]  How can we measure
distance from the simplex $\Delta(\Gamma)$
spanned by a collection $\Gamma$ of disjoint
simple closed curves? 

Define $K(\Gamma)$ to be the set of measured laminations such
that for every $\gamma \in \Gamma$,
$i(\lambda, \gamma)$ = 0.  An element $[\mu] \in \PL_0(S)$
is in $\Delta(\Gamma)$ if and only if $i(\lambda, \mu) = 0$ for
all $\lambda \in K(\Gamma)$.

We can use a normalized version of this to define
neighborhoods of $\Delta(\Gamma)$.  That is, for any
lamination $\lambda \in K(\Gamma)$, the function
$n_\lambda (\mu) = i(\lambda, \mu) / \length_S (\mu)$ depends only on
the class $[\mu] \in \PL_0(S)$.  A neighborhood
for $\Delta(\Gamma)$ is defined by picking a finite
collection of such $\lambda$, and requiring that the
functions $n_\lambda$ have value less than $\epsilon$.

The proof works by induction.  Let $\mu$ be as hypothesized.
We start with the last $i$ such that $i(\alpha_i, \mu) > 0$.
The image of $\mu$ by a large power of this Dehn twist is
near $\alpha_i$, by \ref{big twists}.

Define $\Gamma_j$ to be the set of $\alpha_k$ where $j \le k$,
and for no $j \le l < k$ is $i(\alpha_l, \alpha_k) > 0$.
Assume by induction on $i-j$ that the image $\mu_j$ of $\mu$ by
the $j$ through $n$ terms of the composition is near
$\Delta(\Gamma_j)$.  Consider the result of
of applying the $q-1$st power of the Dehn twist about
$\alpha_{j-1}$ to $\mu_j$.  If the intersection number of
$\mu_j$ with $\alpha_{j-1}$ is not too small compared with
the length of $\mu_j$, then according to
\ref{big twists} the image will be close to $\alpha_{j-1}$.

Otherwise, $\mu_j$ is already near the simplex $\Gamma_{j-1}$.
The length of a lamination can be estimated to within a constant
factor as the sum of the intersection numbers with any
finite collection of simple closed curves $C$ such that $\ls(C) = S$.
Choose $C$ so that it has at least one element which
intersects
any component of $\Gamma_{j-1}$ but none of the others.
Then, if $\mu_j$ is close to the subsimplex of $\Gamma_{j-1}$
opposite $\alpha_j$, we see that Dehn twists about
$\alpha_j$ cannot diminish the length of $\mu_j$ beyond
a bounded factor, since these twists do not affect
certain of the intersection numbers which at the beginning
contribute a significant fraction of the estimate of
the length of $\mu_j$.  Dehn twists about $\alpha_j$ do
not affect the intersection numbers with $\lambda \in K(\Gamma)$,
so they cannot increase the functions $n_{\lambda}$ beyond
a bounded factor.  We conclude that $\mu_{j-1}$
is near $\Delta(\Gamma_{j-1}$.

The proof of the assertion for a metric $g$ is similar,
but slightly simpler, since we can start the induction
at the last term in the composition.  The initial
step follows form \ref{big twists}.
One can measure a neighborhood of a lamination in
compactified Teichm\"uller space in the
analogous way.  The easiest way to deal with the
normalization is to choose a set of curves whose
least subsurface is all of $S$, and normalize a metric
by multiplying with a constant which
makes the sum of their lengths equal to $1$.

A more elegant way is to extend
the function $\length_g(\lambda)$ to a function
$\length_g(h)$, where $h \in \teich(S)$, with the aid of
``random geodesic'' for $g$ and for $h$.  The length
of the $h$-random geodesic in the metric $g$ is the desired
function.  This quantity is also
the intersection number between the random geodesics
in the $h$ metric and the $g$ metric, so it is
symmetric in $g$ and $h$.
See \cite{Th5} for details.

The subsequent steps are identical to the
case for a lamination, since the image
of $g$ is near a lamination.
\end{proof}

The first phase of the proof of Theorem \ref{drill holes}
will be to show that
there is a subsequence such that $\rho_k$ converges.

Before proceeding in general,
it is worth describing how this
phase of the proof goes in some relatively
simple cases.

First consider the case that $S$ is a punctured torus.  Then
$\prec$ is actually a linear ordering.  Consider the surface
at height $n  + 0.5$ --- label this cut $K_n$.
By \ref{compose big twists}, the
conformal structure of the
top component of the domain of discontinuity for $\Gamma(q_j, K_n)$
converges to $\alpha(n+1)$, provided we choose $q_j$ appropriately,
and the conformal structure on
the bottom component
converges to $\alpha(n)$.  The double limit theorem applies,
yielding a subsequence such that the representations restricted
to $S_{K_n}$ converge.
By a diagonal process, we reach a subsequence so that these
groups converge, for all $n$.

For any two consective integers, the image of the fundamental
group of the two surfaces in $M_C$ intersect on the fundamental group
of a thrice-punctured torus.  This implies that the fundamental
group of their amalgamation also converges (since the intersection
is non-elementary).  Inductively,
it follows that the representations of the entire fundamental group of
$M_C$ converges.

Another special case which is fairly easy to handle is the case of an
arbitrary surface $S$ when the
the sequence of curves $C$
is either finite or semi-infinite.  Suppose, for example,
that the index set is the positive integers.  Let $m$ be the
set of minimal elements of $\prec$.  They form a system
of disjoint curves on $S$.  Let $R_1 \subset M_C$ be a surface which is
above all curves in $m$, and below all others.  For appropriate
power functions $q$, the conformal structure at
$+\infty$ will appear as some convex combination of curves
in a collection  $\Gamma$
which are minimal in $C - m$
from the point of view of $S_1$.  Let $m_1 \subset m$ be the subset
of $m$ which elements of $\Gamma$ intersect.  Let $S_1$ be a surface
which is above the $m_i$, but below all other curves in $C$.  
For each element $\beta$ of $m_1$, there is a simple closed curve on $S$ which
intersects $\beta$ but none of the other curves in $m_1$.  From
the point of view of $S_1$, this curve appears as a lamination
close to $\beta$.  Forming some convex
combination of such curves on $S$, we see that there is a lamination
which has modest length
on $S$, and when transformed to $S_1$ is near the sum of the elements
of $m_1$.  An application of
\ref{convergence on non-binding laminations} shows that the subsurface filled up by $\Gamma \union m$
has bounded geometry.  Since the geometry of $S_1 - m_1$ also
remains bounded, the geometry of $S_1$ remains bounded.

This argument can be repeated inductively, using $S_1$ in place of $S$,
to prove that in this circumstance, all of $\fund{M_C}$ converges.

This argment could be extended to the case that $C$ is doubly infinite,
provided
one could find at least one surface whose geometry remained bounded.
For instance, if there is at least one cut in $C$ such that any maximal curve
below the cut together with any minimal curve above the cut
bind the entire surface, that serves as a good starting point.
Examples where this happens are easy to arrange.  Unfortunately,
it is also easy to arrange examples where it is difficult to get
started.  The trouble is that laminations or hyperbolic structures
at $\pm \infty$, when transformed to a surface buried in the middle
of $M_C$, tend to be concentrated near a single curve, rather
than to have a positive coefficient for laminations near
each of the possibile
curves.  Thus, it is hard to control things so that a
lamination which is short at $+\infty$ is
guaranteed to have a positive intersection number with a lamination
which is short at $-\infty$.

Note also that this special case feeds into the remark at the
beginning of the proof of \ref{drill holes}. It is not hard to show
that if $C$ is finite, but large enough to bind the surface, then
$AH(M_C)$ is compact, with an argument similar to the one we have just seen.
Therefore, as $C$ increases through the finite subsets of a bi-infinite
example, the representations remain bounded, so we could arrange
for them to converge to a limit.  We could then diagonalize,
to prove the theorem; this method, however, would lose any information
concerning rates of convergence.

\medskip

So that we can gracefully handle the general case,
we will make a more careful study of laminations and
combinations of laminations which
can have moderate length on hyperbolic structures on $S$.

First, let $\alpha$ be a simple closed curve on $S$.
We will
study laminations which have a small intersection number with
$\alpha$, or in other words, are near $K(\alpha)$.  The
condition $i(\alpha, \lambda) = 0$ is a codimension one condition
on $\lambda$ --- $K(\alpha)$ has dimension one less than
the dimension of $ML_0(S)$.  The concrete
application we will make of this is that if $\mu \in K(\alpha)$
has a positive weight on $\alpha$ itself, then a small neighborhood
of $\mu$ is subdivided into three types of laminations:
\begin{description}
\item[(1)] laminations in $K(\alpha)$, having a positive weight on $\alpha$
\item[(2)] laminations $\nu$ {\it rightward} of $K(\alpha)$, such that 
$i(\nu, \alpha) > 0$, and leaves of $\nu$ spiral toward the right
as they approach $\alpha$, and
\item[(3)] laminations $\nu$ {\it leftward} of $K(\alpha)$, such that 
$i(\nu, \alpha) > 0$, and leaves of $\nu$ spiral toward the left
as they approach $\alpha$.
\end{description}

This subdivision can be recognized with the aid of a short arc $a$
transverse to $\alpha$.  The return map for leaves of laminations
$\nu$ near $\mu$ gives the information:  if leaves return
a bit to the left of where they started, then
one is in the leftward portion of the neighborhood. 
Note that the definitions depend only
on the orientation of $S$, not on an orientation for $\alpha$.

\begin{proposition}[Centrist constriction]\label{centrist constriction}
Let $\alpha$ be a simple closed geodesic on $S$, and let $\mu, \mu' \in
K(\alpha)$ be elements having weight $1$ on $\alpha$. There are neighborhoods
$U$  and $U'$ of $\mu$ and $\mu'$ such that for any hyperbolic structure $g$ on
$S$, if there is a lamination $\nu$ in $U$ rightward of $K(\alpha)$ and a
lamination $\nu'$ in $U$ leftward of $K(\alpha)$ such that $\length_g (\nu) <
A$ and $\length_g(\nu') < A$, then there is another lamination $\lambda \in
K(\alpha)$ having weight $1$ on $\alpha$, such that
$$ \length_g(\lambda) < 1.1 A$$.
\end{proposition}
\begin{proof}[Proof of \ref{centrist constriction}]
The lamination $\lambda$ will come from the ``complex geodesic'' determined by
$\nu$ and $\nu'$.  The laminations $\mu$ and $\mu'$ are really only used to
control the slope of leaves of $\nu$ and $\nu'$ near $\alpha$. We want the
slopes to be small. If $\nu$ and $\nu'$ have any leaves in common, throw out
the union of all such leaves, to obtain new laminations whose leaves are still
nearly tangent to $\alpha$.  Also throw out any leaves which do not meet the
other lamination.

Once the leaves of $\nu$ and $\nu'$ are transverse, there is a construction for
laminations which are weighted combinations of $\nu$ and $\nu'$.  This can be
done with the help of weighted train tracks. Construct a foliation $\mathcal F$ in
a neighborhood of $\nu \union \nu'$ which is transverse to both, using the
condition that if local coordinates are chosen so that the leaves of $\nu$ are
horizontal and leaves of $\nu'$ are vertical, the leaves of $\mathcal F$ run from
the lower left to the upper right.  In a neighborhood of $\alpha$, $\mathcal F$ can
be taken perpendicular to $\alpha$. If the support neighborhood for $\mathcal F$ is
trimmed down far enough, then all its leaves are intervals and bounded in
length. The intervals fall into a finite number of isotopy classes within the
support neighborhood.  The quotient space of the neighborhood by the leaves of
$\mathcal F$ is a train track --- a good quality train track (in that the
complementary regions have the right types) if the neighborhood is trimmed
carefully enough so that each leaf of $\mathcal F$ meets a leaf of $\nu$ or of
$\nu'$.  Each of $\nu$ and $\nu'$ are carried on the resulting train track with
positive weights.   The lamination $\lambda(s,t)$ is obtained by taking the
$(s, t)$ convex combination of the weights for $\nu$ and the weights for $\nu'$.

Clearly, the length of $\lambda(s,t)$ does not exceed $s$ times the length of
$\nu$ plus $t$ times the length of $\nu'$. Choose the ratio of $s$ and $t$ so
that the leftward shifting and the rightward shifting of the return map to $a$
exactly cancel.  The resulting lamination has a positive weight on $\alpha$. 
Its weight can be estimated by looking at the total transverse measure of
intersection with $\alpha$. This weight is approximately $s + t$, but the
formula is not precise.  Multiplying by a constant near $1$, we obtain the
lamination $\lambda$ as claimed.
\end{proof}

\begin{corollary}[Curves twisted tight]\label{curves twisted tight}
Let $C$ be a family of curves in $S \cross \reals$, as above, and $q$ as above
an appropriate choice of Dehn surgeries. Let $K \prec K'$ be cuts which differ
by only one curve $\gamma = \alpha \cross t$.  Suppose that $\nu$ is a
lamination which has moderate length on the top domain of discontinuity from
the point of view of $S_K$, and $\nu'$ is a lamination which has moderate
length on the bottom domain of discontinuity from the point of view of
$S_{K'}$.  If $\nu$ and $\nu'$ are close to laminations in $K(\alpha)$ with
definite components of $\alpha$, then the length of $\alpha$ in the quotient
three-manifolds of $\Gamma(q, K)$ and $\Gamma(q, K')$ is moderate.
\end{corollary}
\begin{proof}[Proof of \ref{curves twisted tight}]
The lamination $\nu$ has the propoerty that a certain composition $\phi$ of
high powers of Dehn structures about a sequence of curves, beginning with
$\alpha$, sends it to a bounded lamination. Let us use the convention that the
$\tau_\alpha$ twists toward the left, so that a high positiive power of
$\tau_\alpha$ sends most laminations to leftward laminations near $K(\alpha)$.
If the power of $\tau_\alpha$ which occurs as the first term in $\phi$ is
positive, then $\nu$ must be rightward of $K(\alpha)$, by
\ref{compose big twists} and its proof --- each step in the proof either
maintains the lamination roughly unchanged, or drastically increases its
length.  The only chance for $\nu$ to transform to something of moderate length
is for it to reduce in length at the first step, which will only happen if it
is a rightward lamination.

Similarly, if the power of $\tau_\alpha$ which occurs as the first term in
$\phi$ is positive, then the power of $\tau_\alpha$ which occurs as the first
term in the analogous diffeomorphism $\phi'$ which sends $\nu'$ to a lamination
with moderate length is negative, and $\nu'$ must be a leftward lamination.

The two laminations $\nu$ and $\nu'$ are on different surfaces, differing by a
change of parametrization of the form $\tau_\alpha^N$. There is a compromise
parametrization, obtained by applying the $[N / 2 ]$ power of $\tau_\alpha$ to
one of the two.  When $\nu$ and $\nu'$ are transformed to the compromise
surface, they are near laminations with roughly half the size of their
approximate components of $\alpha$, and the images are still on opposite sides
of $K(\alpha)$.  If we represent this surface by an efficient pleated surface,
we obtain both leftward and rightward laminations on the same surface, so we
can apply \ref{centrist constriction} to conclude that on this surface,
$\alpha$ itself has moderate length.

\end{proof}

\begin{proposition}[Minor detours]\label{minor detours}
Let $h$ be a hyperbolic structure fixed on $S$ for reference,
and let $\alpha$ be a simple closed geodesic on $S$.

For every $\epsilon$ there is a $\delta$ such that if $g$ is a hyperbolic
structure on $S$ such that
$$\length_g(\alpha) < A \length_h(\alpha),$$
and if $\mu \in \ML_0(S)$ is a measured lamination normalized to have length 1
such that
$$i(\mu, \alpha )< \delta$$
and
$$\length_g(\mu) < B,$$
then $\mu$ is near a lamination $\mu' \in \ML_0(S)$ of length $1$ such that
$$\length_g(\mu') < (1+\epsilon)(A+B) $$
and
$$i(\mu', \alpha) = 0.$$
\end{proposition}

\begin{proof}[Proof of \ref{minor detours}]
We will modify $\mu$ by cutting the strands of $\mu$ which intersect $\alpha$,
and joining the resulting ends by paths which run near intervals of $\alpha$.

We begin by performing the positive or negative power $P$ of $\tau_\alpha$
which reduces the length of $\mu$ as much as possible.  At the very end we will
compensate for this unwinding around $\alpha$ by adding a component $i(\alpha,
\mu) |P| \alpha$ to obtain the final result.

Let $\mu_1$ be the result of this first unwinding step.  If $\mu_1$ has a very
small length, we can proceed immediately to the final step; a multiple of
$\alpha$ is then sufficiently close to $\mu$.

Otherwise, proceed as follows. Some components of $\mu_1$ may intersect
$\alpha$, and others may not. Since every component of a measured lamination is
minimal, every leaf of the components which do intersect $\alpha$ intersects in
a bi-infinite sequence. The leaves when cut by $\alpha$ group into a finite
number of bands of isotopic and parallel arcs, where the number is bounded by
the genus of the surface.  The angles at which the leaves meet $\alpha$ are
bounded away from zero (depending on the geometry of $h$).

The lengths of arcs in a given band differ by less than $2 \length_h
(\alpha)$.  If $\delta$ is small compared to the length of $\mu_1$, then the
average length in $h$ for bands constituting most of the total transverse
measure of $\mu_1$ along $\alpha$ is large. We assume that $\delta$ is
sufficiently small, and we begin the process of forming $\mu'$ by cutting all
the leaves of $\mu_1$ along $\alpha$, retracting their ends a bit so that they
end on the boundary of a thin annular neighborhood of $\alpha$, and throwing
all bands whose length in $h$ is less than some fairly large constant.

We will describe the construction of $\mu'$ by constructing a weighted train
track which carries it. Choose one representative arc from each band, and mark
it with the total transverse measure for the band. Now isotope all the arcs so
as to minimize their length.  The result of the isotopy will be a family of
arcs which are orthogonal to $\alpha$. The isotopy moves the endpoints a
bounded distance, since the angles were bounded at the beginning. Consequently,
most of the length of each of the arcs moves a rather short distance.

The endpoints of the arcs divide each side of the annulus into a bounded number
of intervals, so the longest interval has length bounded above zero.  Join the
ends of the arcs which bound the longest interval by an arc inside the annulus,
and assign it a weight which is the minimum of the weights for the two arcs. 
If two weights were different, then extend the arc with the larger weight,
assigning the difference to the extension.  Now make an isotopy to retract the
added portion out of the annulus, until only the one possible endpoint touches
the boundary of the annulus.  Continue inductively.  At the last step, there
may be only one arc touching a boundary component of the annulus; in that case,
add a loop which goes completely around the boundary of the annulus, with half
the weight of that arc.

Each trajectory on the weighted train track we have constructed is homotopic to
a polygonal path, with right angle bends, with side lengths (in $h$)
alternately bounded below by a constant greater than zero, and bounded below by
a fairly large constant.  Such a path is homotopic to a geodesic which is close
to the original polygon, along most of the lengths of the long segments.  This
means that the weighted train track we have constructed determines a lamination
$\mu_2$ which is near $\mu_1$, and such that $\length_h(\mu_2)$ is close to
$\length_h(\mu_1)$.

On the other hand, we can represent the train track for $\mu_2$ in $g$ by using
leaves of $\mu$, joining them by short arcs near $\alpha$.  The total length of
the resulting weighted train track is less than $B + i(\alpha, \mu) A 
\length_h(\alpha)$.

Now add $i(\alpha, \mu)P\alpha$ to $\mu_2$, to obtain $\mu'$.  The
resulting lamination is close to $\mu$, and the length of its
image in $h$ satisfies the inequality asserted in the proposition.
\end{proof}

We return now to the proof of \ref{drill holes}. We will focus on an
arbitrary element $\gamma \in C$, and show that its length remains bounded. 

Start with the cut $G$ whose lower part is the set of all curves $\beta \prec
\gamma$ together with $\gamma$ itself. It follows from the hypotheses of the
theorem that this in fact defines a cut.  The lower half of $G$ has a unique
maximal element, namely $\gamma$, so the lower conformal structure from the
point of view of $S_G$ is near $\gamma$, by \ref{compose big twists}.  Let
$\mu_0$ a lamination near $\gamma$ which has very small length in this
conformal structure, guaranteed by Theorem
\ref{laminations compactify Teichmuller space}.

We will construct an increasing sequence $K_i$ of cuts for which $\gamma$ is a
maximal element in the lower part, and such that there is a lamination $\mu_i$
on $K_i$ which is  either short in the bottom conformal structure, or at least
has bounded length in $N$, and which has a substantial approximate component of
$\gamma$.  We will continue until we reach a cut where the upper conformal
structure is near a lamination having positive intersection number with
$\gamma$. There is a greatest cut $H$ for which $\gamma$ is a maximal element
in the lower part: its upper part consists of the curves $\beta$ such that
$\gamma \prec \beta$. If we ever reach $H$, we will be done, since all minimal
elements of the upper part of $H$ intersect $\gamma$.

Set $K_0 = G$.  When $K_i$ has been defined, let $\lambda_i$ be a lamination
whose length is near 0 in the upper conformal structure from the point of view
of $R = S_{K_i}$.  Let $\beta_i$ be a linear combination of simple closed
curves on $R$ isotopic to minimal elements above $R$ which closely approximates
$\lambda_i$, and let $\delta_i$ be the curve with the largest weight.

If $i(\delta_i, \gamma) > 0$, then stop, construct an efficient pleated surface
representing $R$, and apply \ref{laminations bind subsurface}, to
conclude that $\gamma$ has bounded length in $N$.

Otherwise, define $K_{i+1}$ by moving $\delta_i$ from the upper part to the
lower part, and abbreviate $Q = S_{K_{i+1}}$ Let $\nu$ be $\mu_i$ pulled back
in the homotopy class of $R$ rather than $Q$. There are two cases.

First, it may happen that $\nu$ has a substantial approximate component of
$delta_i$.  In that case, observe  that $\nu$ and $\delta_i$ must be on
opposite sides of $K(\delta_i)$, as in proposition
\ref{curves twisted tight}. Therefore, $\delta_i$ has bounded length in
$N$.  Let $R$ be represented efficiently by some pleated surface, so that
$\mu_i$ has bounded length on $R$.  Since $\delta_i$ also has bounded length,
we can  apply \ref{minor detours}, to find another lamination $\mu'$ near
$\mu_i$ which has bounded length on this pleated surface.  Define $\mu_{i+1}$
to be $\mu$ transported to $Q$.  We may as well assume that the component of
$\delta_i$ in $\mu_{i+1}$ is zero, since it is negligibly small anyway.

The other case is that $\nu$ does not have a substantial component of
$\delta_i$.  In that case, we can take $\mu_{i+1} = \nu$.  The domain length
(on $S$) of $\mu_i$ and of $\mu_{i+1}$ is approximately the same, and in the
range, the laminations are exactly the same.

In the end, we must eventually arrive at a point where $i(\delta_i, \gamma) >
0$.  This shows that the length of $\gamma$ is bounded.

It is worth commenting here on the constants involved in this induction,
because it may be useful someday for the development of more precise analysis
of the quasi-isometric geometry of surface groups.

How does the ratio $B$ of length of $\mu_i$ in $N$ to that in $S$ progress? 
The ratio starts out near $0$.

If the power of the Dehn twist about $\gamma$ is reasonably large, the
component of $delta_i$ in $\mu_i$ is negligible in the current situation, so
the constant $(1+\epsilon)(A+B)$ could be replaced by $(1+\epsilon)B$ in a step
of the first type.

We lose an additive constant $a(\lambda, \mu_i)$
when we apply \ref{efficiency of pleated surfaces}.
However, for the purposes of this proof, we could suppose that $\gamma$ does
not have a short length (even though it does).  In that case, we could take
$\lambda$ always to have a closed leaf isotopic to $\gamma$ --- this would mean
that $a(\lambda, \mu_0)$ at least is quite small.  A similar trick works at
successive stages --- if $\mu_i$ has approximate components on curves which are
not short, choose $\lambda$ to contain closed curves in those homotopy classes,
otherwise use \ref{minor detours} to eliminate those approximate
components.

\medskip

Continuing with the proof of \ref{drill holes}, we take a subsequence of
the representations of $\fund {M_C}$ such that restricted to any surface $S_K$,
the subsequence is maximally convergent.  What are the subsurfaces of maximal
convergence?  An argument similar to the previous argument, for the case that
the sequence of curve is finite or semi-infinite, shows that representations
converge on the entire surface $S_K$, as follows. The subsurface $R$ of
convergence for $S_K$ must contain curves in the isotopy classes of at least
the set $m$ of all minimal elements above $S_K$ and the set $M$ of all maximal
elements below. Let $\delta$ be a minimal element greater than the set $m$, and
let $m_1 \subset m$ be those elements of $m$ which are less than $\delta$. 
Consider the surface $S_{K'}$, where $K'$ is obtained from $K$ by moving $m_1$
to the lower part.  The subsurface $R'$ of convergence for $S_{K'}$ contains
$R$ cut by $m_1$, and also contains $\delta$. We can continue upward in this
way, until at some point we arrive at a surface whose subsurface of convergence
is the entire surface.  As we have already seen, this implies that the entire
fundamental group of $M_C$ converges.

\bigskip

The second phase of the proof of \ref{drill holes} is to prove that the
geometric limit of the manifolds is, in fact, homeomorphic to $M_C$.

We claim first that the limit representation of $\fund{M_C}$ is faithful.  It
is a general fact that if the limit a sequence of representations of a group
$G$ is non-faithful, and if $G$ is not almost abelian, then any element $g$
which is represented trivially in the limit must be represented trivially for
all but a finite number of representations in the sequence.  We know exactly
the kernel of the $i$th representation: it is the same as the kernel of the
homomorphism $\fund {M_C}$ to $\fund S$.  It is easy to verify (e.g., by
estimating the hyperbolic geodesic on $S$ representing the image of a loop in
$M_C$) that the shortest word in the kernel has length which goes to infinity
with $i$, so the limit representation indeed is faithful.

Let $P$ be the geometric limit manifold, and let $N$ be the quotient of $\hy^3$
by the limiting representation of $\fund {M_C}$. By general principles, there
is a covering map $N \arrow P$. We must show that the covering is trivial, and
that $N$ is indeed homeomorphic to $M_C$.

The covering cannot be finite to one: this can be seen with the aid of the
unique divisibility property for $\fund S$ and $\fund {M_C}$, that if any
element is expressible as a $k$th power, it is expressible in only one way as a
$k$th power. Furthermore, if an infinite number of the images of $\beta \in
\fund {M_C}$ in $\fund S$ are divisible by $k$, then $\beta$ itself is
divisible by $k$. Therefore, if $\alpha \in \fund P$ has $k$th power in $\fund
N$, it follows that $\alpha \in \fund N$.

The other possibility might be that the covering is infinite to one.  This
possibility we will rule out by a geometric argument. First, by
\refin{Ch. 8 and 9}{GT3M}, the injectivity radius is bounded above
throughout $N$. This implies that in any infinite-sheeted covering projection,
cusps are identified only with cusps.  The projection of each cusp to $P$ is
one-to-one, since it is finite-to-one.

Let $\Gamma \subset C$ be an independent set of maximal size. Construct a
pleated surface $P$, obtained by pinching all elements of $\Gamma$ to their
cusps. 

Consider first the case that $C$ is bi-infinite. The space $B$ of $2$-cycles in
$N$, up to chains of bounded volume, is isomorphic to $\integers$, and any
surface $S_K$ represents a generator.

If the covering $N \arrow P$ is infinite to one, then each cusp is a member of
an infinite sequence of cusps identified to one cusp in $P$. There are pleated
surfaces representing a generator of $B$ passing through each of the cusps, by
\refin{Ch. 8 and 9}{GT3M}.  If some of these surfaces have $\epsilon$-thin
parts which are isotopic to cusps of $N$, we can modify them by actually
pinching such short curves to their cusps.  But all surfaces of this type are
contained in a subset of $P$ of bounded volume and finitely generated
fundamental group, again by  \refin{Ch. 8 and 9}{GT3M}, since $N$ and
hence $P$ has no $\integers$-cusps.  (Any two points on a pleated surface are
connected by a path whose total length, excluding its intersection with the
thin part of the surface, is {\it a priori} bounded.  The set of points in $P$
accessible by such paths has bounded volume and finitely generated fundamental
group). It follows that some of the pleated surfaces are mapped to isotopic
pleated surfaces in $P$, and that the total algebraic volume of the isotopy is
bounded.  Combining this isotopy with the image of a chain in $N$ whose
boundary is the difference of the cycles represented by the two pleated
surfaces, we obtain a 3-cycle of finite but non-zero volume in $P$ --- $P$ has
finite volume.

However, the geometric limit of hyperbolic manifolds of infinite volume never
has finite volume --- a limit manifold is homeomorphic to the approximating
manifolds in its thick part, which implies that any approximants also would
have finite volume.

Therefore, the covering $N\arrow P$ is trivial, and the geometric limit is $N
\homeomorphic M_C $, as claimed.

If the sequence is not bi-infinite, we can argue with the domain of
discontinuity.There is say a lowest surface or a highest surface $Q$.  In the
sequence of representations for $\fund Q$, the conformal structure on the
bottom domain of discontinuity is constant.  In the limit, this bottom domain
of discontinuity is still present.  It embeds in the domain of discontinuity
for the geometric limit.  Thus there is a covering map $D_{\fund N} /\fund N
\arrow D_{\fund P}/\fund P$.  Since both domain and range are hyperbolic
surfaces of finite area, such a covering map can only be finit-to-one, so as we
have already seen, it must be one-to-one.

We have a homotopy equivalence $f\colon M_C \arrow N \homeomorphic P$. Why is
it homotopic to a homeomorphism? First, we can easily homotope until $f$ is a
homeomorphism near the cusps. For each surface cut $S_K$, the fundamental group
is represented as a free product amalgamated along $\fund S_K$, so by standard
three-manifold topology, we can find an incompressible surface whose fundmental
group is contained in $\fund S_K$ and which separates $N$ into two pieces with
the expected fundamental groups, and which meets the cusps of $N$ in the
expected way. The only possibility, for homological reasons, is that the
surface has the full fundamental group of $S_K$, and so it is homeomorphic to
$S_K$.  We can continue, by the standard technique, to conclude that $f$ is
homotopic to a homeomorphism $M_C \homeomorphic N$.
\end{proof}

\bibliographystyle{halpha}
\bibliography{hype}
\end{document}